\documentclass[graybox]{article}

\usepackage{type1cm}        
\usepackage{graphicx}        

\usepackage{newtxtext}       %
\usepackage{newtxmath}       

\usepackage[final]{optional}
\usepackage{amsfonts, amsmath, wasysym}
\usepackage{tikz}
\usepackage{amssymb, 
paralist}
\usepackage{latexsym}
\usepackage{url}

\usepackage{pgfplots}

\newtheorem{lemma}{Lemma}
\newtheorem{theorem}{Theorem}

\newtheorem{definition}{Definition}

\newtheorem{conjecture}{Conjecture}

\newtheorem{example}{Example}
\newtheorem{remark}{Remark}

\newcommand{\s}{\ensuremath{{\cal S}}}

\newcommand{\cc}{\ensuremath{{\cal C}}}
\newcommand{\cd}{\ensuremath{{\cal D}}}

\newcommand{\calr}{\ensuremath{{\cal R}}}

\newcommand{\pl}[2]{{\frac{\partial #1}{\partial #2}}}

\newcommand{\pll}[2]{{\frac{\partial^2 #1}{\partial #2^2}}}

\newcommand{\ti}{\tilde}

\newcommand{\al}{\alpha}
\newcommand{\be}{\beta}
\newcommand{\ga}{\gamma}

\newcommand{\de}{\delta}
\newcommand{\om}{\omega}
\newcommand{\Om}{\Omega}

\newcommand{\la}{\lambda}

\newcommand{\si}{\sigma}

\renewcommand{\th}{\theta}

\newcommand{\vph}{\varphi}
\newcommand{\ep}{\varepsilon}

\newcommand{\R}{\ensuremath{{\mathbb R}}}
\newcommand{\N}{\ensuremath{{\mathbb N}}}
\newcommand{\B}{\ensuremath{{\mathbb B}}}

\newcommand{\downto}{\downarrow}
\newcommand{\upto}{\uparrow}

\newcommand{\lap}{\Delta}

\newcommand{\grad}{\nabla}

\DeclareMathOperator{\Vol}{Vol}
\DeclareMathOperator{\VolB}{VolB}

\DeclareMathOperator{\inj}{inj}

\def\blbox{\quad \vrule height7.5pt width4.17pt depth0pt}

\newcommand{\beq}{\begin{equation}}
\newcommand{\beql}[1]{\begin{equation}\label{#1}}
\newcommand{\eeq}{\end{equation}}
\newcommand{\beqa}{\begin{equation}\begin{aligned}}
\newcommand{\eeqa}{\end{aligned}\end{equation}}
\newcommand{\brmk}{\begin{remark}}
\newcommand{\ermk}{\end{remark}}
\newcommand{\partref}[1]{\hbox{(\csname @roman\endcsname{\ref{#1}})}}
\newcommand{\half}{\frac{1}{2}}

\newcommand{\cmt}[1]{\opt{draft}{\textcolor[rgb]{0.5,0,0}{
$\LHD$ #1 $\RHD$\marginpar{\blbox}}}}

\newcommand{\tr}{{\rm tr}}

\newcommand{\Rm}{{\mathrm{Rm}}}
\newcommand{\Ric}{{\mathrm{Ric}}}

\usepackage{soul}

\begin{document}

\title{Ricci flow and Ricci Limit Spaces}
\author{Peter M. Topping}
\date{6 May 2020}
%
%

\maketitle


%
%



\section*{Preface}
In these lecture notes I hope to describe some of the developments in Ricci flow theory from the past decade or so in a manner accessible to starting graduate students. I will be most concerned with the understanding of Ricci flows that are permitted to have unbounded curvature in the sense that  the curvature can blow up as we wander off to spatial infinity and/or as we decrease time to some singular time. This is an area that has seen substantial recent progress, and it has long been understood that a sufficiently developed theory will be  useful in settling important open problems in geometry, particularly to the subject of understanding the topology of manifolds that satisfy certain curvature restrictions. In these notes I will outline an application of the theory as developed so far, but instead to the theory of so-called Ricci limit spaces, as developed particularly by Cheeger-Colding around 20 years ago. 

Much of the work I describe is joint work with Miles Simon. Some of the general ideas will be explained in the context of two-dimensional flows -- the case in which these were developed -- which is partly joint work with Gregor Giesen.

These notes have been written for a summer school at Cetraro in 2018 and I would like to thank Matt Gursky and Andrea Malchiodi for the invitation to explain this theory. A subset of the lectures were given at the Max Planck Institute for Mathematics in Bonn in 2017, and I would like to thank Karl-Theodor Sturm for that invitation.

\vskip 0.3cm
\noindent
{\em
Warwick, 21 December 2018
\hfill
Peter Topping
}

\tableofcontents

\section{Introductory material}

\subsection{Ricci flow basics}

Traditionally, the Ricci flow has been posed  as follows. Given a Riemannian manifold $(M,g_0)$, we ask whether we can find a smooth one-parameter family of Riemannian metrics $g(t)$ for $t\in [0,T)$ such that $g(0)=g_0$ and satisfying the PDE
\beql{rf_eq}
\pl{g}{t}=-2\Ric_g,
\eeq
where $\Ric_g$ is the Ricci curvature of $g$.
If you are new to this subject then you should view the right-hand side of \eqref{rf_eq} morally as a type of Laplacian of $g$, and the equation as an essentially parabolic PDE for the metric $g$.
One of the simplest examples that is not completely trivial is the shrinking sphere: if $(M,g_0)$ is the round unit sphere of dimension $n$, then one can easily check that 
\beql{shrink_sphere}
g(t)=\big(1-2(n-1)t\big)g_0
\eeq
is a solution for $t\in [0,\frac{1}{2(n-1)})$.
Another standard example that is important to digest is the dumbbell example: If we connect two round three-spheres by a thin neck $\ep S^2\times [0,1]$ to give a nonround three-sphere, then Ricci flow will deform it by shrinking the thin neck, which will degenerate before much happens to the two original three-spheres. See \cite[Section 1.3.2]{RFnotes} for a more involved discussion and pictures.

\cmt{will need the picture of a neckpinch for Example \ref{ex4}}

The Ricci flow is a little simpler if the evolving manifold is K\"ahler, and the simplest of all is in the lowest possible dimension, i.e. $\dim M=2$, which corresponds to one complex-dimensional \emph{K\"ahler Ricci flow}. 
Then we can write $\Ric_g=Kg$, where $K$ is the Gauss curvature, and the flow can be seen to deform the metric conformally.

The Ricci flow was introduced by Hamilton in 1982 \cite{ham3D} in order to prove the following theorem.\footnote{When we write $\Ric>0$, we mean that $\Ric(X,X)>0$ for every nonzero tangent vector $X$.}

\begin{theorem}[Hamilton \cite{ham3D}.]
If $(M,g_0)$ is a closed, simply connected Riemannian 3-manifold with $\Ric>0$, then it is diffeomorphic to $S^3$.
\end{theorem}
The rough strategy is as follows. First, take $(M,g_0)$ and run the Ricci flow. Second, show that the condition $\Ric>0$ is preserved for each of the subsequent metrics $g(t)$. 
(It is worth remarking that a local form of this preservation of lower Ricci control has recently been proved \cite{ST1} and will be pivotal later on; see Lemma \ref{modifieddouble2}.) 
Finally, show that after appropriate scaling, $(M,g(t))$ converges to a manifold of constant sectional curvature as time increases. This must then be $S^3$.

\cmt{need $c/t$ curvature decay for Lemma \ref{modifieddouble2} above}

For more details of this argument, from a modern point of view, and a more involved introduction to Ricci flow, see \cite{RFnotes}.

This application is a model for many of the most remarkable subsequent applications of Ricci flow. That is, if we assume that a manifold satisfies a certain curvature condition, then often  we can use Ricci flow to deform it smoothly to something whose topology we can identify (or to a \emph{collection} of objects whose topology we can identify). The famous work of Perelman goes one step further, in some sense, by carrying out that programme in three dimensions without any curvature hypothesis whatsoever, thus proving the Poincar\'e conjecture and Thurston's geometrisation conjecture as originally envisaged by Hamilton and Yau.

In contrast, in these lectures we will use Ricci flow for a completely different task.
We will have to pose Ricci flow in a different way, specifying initial data that is a lot rougher than a Riemannian metric, and as a result we will get a refined description of so-called Ricci limit spaces as we explain later. However, the developments in the Ricci flow theory that permit these new applications seem also likely to lead to 
new applications of the traditional type, i.e. assume a curvature condition and deduce a topological conclusion.

\subsection{Traditional theory - closed or bounded curvature}

The main difference between the theory we will see in these lectures and the traditional theory of Ricci flow is that most traditional applications concern \emph{closed} manifolds, whereas we are interested in general noncompact manifolds. In the closed case, Hamilton proved existence of a unique solution to Ricci flow, as posed above. In the non-closed case, things become mysterious and involved, unless one reduces essentially to a situation like the closed case by imposing a constraint of \emph{bounded curvature}.

\begin{theorem}
[{Shi \cite{shi}, Hamilton \cite{ham3D}.}]
Suppose $(M,g_0)$ is a complete {\bf bounded curvature} Riemannian manifold of arbitrary dimension -- let's say $|\Rm|_{g_0}\leq K$ for some $K>0$. Then there exists a Ricci flow $g(t)$ on $M$ for $t\in [0,\frac{1}{16K}]$ with $g(0)=g_0$ such that $|\Rm|_{g(t)}\leq 2K$ for each $t$.
\label{shi_thm}
\end{theorem}
Here we are writing $\Rm$ for the full curvature tensor. For fixed dimension, controlling $|\Rm|$ is a little like controlling the magnitude of the  sectional curvature.

\cmt{But for a round sphere, $|\Rm|$ divided by the sectional curvature is growing in $n$!}

Because the Ricci flow of Theorem \ref{shi_thm}  has bounded curvature, it will automatically inherit the completeness of $g_0$, i.e. at each time $t$, the metric $g(t)$ will be complete. This is because the lengths of paths can only grow or shrink exponentially:

\begin{lemma}
Given a Ricci flow $g(t)$, $t\in [0,T]$, with
$-M_1\leq \Ric_{g(t)}\leq M_2$,
we have
$$e^{-2M_2t}g(0)\leq g(t) \leq e^{2M_1t}g(0),$$
for each $t$.
\end{lemma}
Here, when we write $\Ric_g \leq K$, say, we mean that 
$\Ric_g\leq K g$ as bilinear forms.

Thus, distances can only expand by a factor at most $e^{M_1t}$, and shrink by a factor no more extreme than $e^{-M_2t}$ after a time $t$.
This lemma is a simple computation. For example, for a fixed tangent vector $X$, we see that 
$$\pl{}{t}\ln g(X,X)\leq 2M_1.$$
We also remark that, as proved by Chen-Zhu \cite{chenzhu} and simplified by Kotschwar \cite{kotschwar}, the flow of Theorem \ref{shi_thm} is unique in the sense that two smooth bounded-curvature Ricci flows starting at the same smooth, complete, bounded curvature initial metric must agree. We emphasise that by \emph{bounded-curvature} Ricci flow we 
mean that there is a $t$-independent curvature bound, not that the curvature of each $g(t)$ is bounded, otherwise the uniqueness fails (see, e.g., Example \ref{ex1} below).

\brmk
\label{Shi_length}
For how long does Shi's flow live? The theorem as we gave it tells us that we can at least flow for a time $1/(16K)$, which is, by chance, independent of the dimension of $M$. But the theorem can be iterated to give a flow that lives until the curvature blows up, and because the theorem provides a curvature bound at the final time, the iterated flow must live for at least 
a time $1/(16K)+1/(32K)+1/(64K)+\cdots=1/(8K)$. None of these explicit times are claimed to be optimal. Traditionally one called a flow whose curvature blew up at a finite end time a \emph{maximal} flow. However, it was shown in \cite{CR_W} and \cite{GT3} that one can have flows in this context whose curvature blows up only at spatial infinity at the end time, and which can be extended beyond to exist for all time (see also
Section \ref{2Dsect}).
\ermk

\cmt{Note that a unit sphere will collapse in a very short time for very high dimension $n$, but will also have a very high $K$...}

\section{Noncompact and unbounded curvature case}

\subsection{Difficulties of the noncompact case; examples to reboot intuition}


Shi's theorem \ref{shi_thm} gives a fine existence statement in the special case that we wish to flow a manifold with bounded curvature. Unfortunately, this is an unreasonable assumption in virtually all applications when we are working on a noncompact manifold.  Thus we would like to consider the problem of starting a Ricci flow with a more general manifold, particularly one of unbounded curvature, or more generally considering flows that themselves have unbounded curvature.
Unfortunately, such flows can be very weird.
They can jump back and forth between being complete and being incomplete. They can exhibit unusual nonuniqueness properties.
They tend to render the ever-useful maximum principle void.
They will  try to jump to a different underlying manifold if given the chance.
Ultimately, they will simply fail to exist in certain cases. 

Let us elaborate on some of these vague statements with some specific examples.

\begin{example}
\label{ex1}
If $(M,g_0)$ is the two-dimensional hyperbolic plane, then there is 
a simple expanding solution $g(t)=(1+2t)g_0$ for $t\in [0,\infty)$ that is the direct analogue of the shrinking sphere that we saw in \eqref{shrink_sphere}. This is the solution that would be given by Shi's theorem \ref{shi_thm}. 
However, more surprisingly, there exist infinitely many other smooth Ricci flows also starting with the hyperbolic plane, and with bounded curvature for $t\in (\ep,\infty)$, for any $\ep>0$. As we will see, these solutions are necessarily  not complete, and thus must have unbounded curvature as $t\downto 0$, although being smooth they trivially have bounded curvature as $t\downto 0$ over any compact subset of the domain.
Intuitively they instantly peel away at the boundary, though not necessarily in a rotationally symmetric way.
We'll construct them in a moment.
\end{example}

\cmt{promising the 2D uniqueness statement above when we say that the solutions are necessarily not complete}

\begin{example}
\label{ex2}
There exists a rotationally symmetric, complete, bounded curvature conformal metric $g_0$ on $\R^2$, and a smooth rotationally symmetric Ricci flow $g(t)$
for $t\in [0,1)$ with $g(0)=g_0$, such that the distance from the origin to infinity (i.e. the length of a radial ray) behaves like $-\log t$ (i.e. it is bounded from above and below by $-\log t$ times appropriate positive constants). In particular, the boundary at infinity is uniformly sucked in making it instantaneously incomplete. The curvature at infinity behaves like $t^{-2}$ (again controlled, up to a constant factor, from above and below). 
\end{example}
Before continuing with our examples, let's construct the two we have already seen. They both have the same building block, which is the concept of a Ricci flow that contracts a hyperbolic cusp. We showed in \cite{revcusp} (though the best way of seeing this now is to use the sharp $L^1$-$L^\infty$ smoothing estimate proved in \cite{TY1}) that if we are given a surface with a hyperbolic cusp, then one way of flowing it is to add a point at infinity and let it contract. 

To see a precise instance of this, consider the punctured two-dimensional disc $D\setminus\{0\}$, equipped with the unique (conformal) complete hyperbolic metric $g_h$. 
That is, 
$$g_h=\frac{1}{r^2\log r^2}(dx^2+dy^2),$$
where $r^2=x^2+y^2$.
This has a hyperbolic cusp at the origin. Similarly to what we have seen before, there is an explicit Ricci flow starting with $g_h$ given by $g(t)=(1+2t)g_h$, and this is the flow that Shi's theorem \ref{shi_thm} would give. However, in \cite{TY2} we constructed a smooth Ricci flow $g(t)$ on the whole disc $D$, for $t\in (0,\infty)$ such that $g(t)\to g_h$ smoothly locally on $D\setminus\{0\}$
as $t\downto 0$. This flow we can restrict to $D\setminus\{0\}$, giving a smooth Ricci flow on $D\setminus\{0\}$ for $t\in [0,\infty)$ such that 
$g(0)=g_h$, but so that $g(t)$ is incomplete at the origin for $t>0$.

To construct Example \ref{ex1}, we merely have to lift this flow on 
$D\setminus\{0\}$ to its universal cover $D$.

Meanwhile, to construct Example \ref{ex2}, we take a slight variant of the $D\setminus\{0\}$ example above by taking a complete rotationally symmetric metric $g_0$ on a punctured 2-sphere, with a hyperbolic cusp.
Again, we can find a rotationally symmetric Ricci flow on the whole 2-sphere for 
$t\in (0,\ep)$, that converges to $g_0$ on the \emph{punctured} sphere as $t\downto 0$, and we can restrict this flow to the punctured sphere, which is conformally $\R^2$. The details and the claimed asymptotics can be found in \cite{revcusp, TY1, TY2}.



In both Examples \ref{ex1} and \ref{ex2} we see incompleteness in the flow itself. In the theory we'll see later it turns out to be important to be able to flow \emph{locally}, in which case even the initial data is incomplete. The following example shows how bad that can be.

\begin{example}
[{Taken from \cite{GT3}.}]
For all $\ep>0$, there exists a Ricci flow $g(t)$ on the unit two-dimensional disc $D$ in the plane, defined for $t\in [0,\ep)$, such that $g(0)$ is the standard flat metric corresponding to the unit disc, but 
$$\inf_{D}K_{g(t)}\to\infty\qquad\text{as }t\upto \ep,$$
and $\Vol_{g(t)}D\to 0$ as $t\upto \ep$.
\label{ex3}
\end{example}
If you've ever seen Perelman's famous Pseudolocality theorem before \cite[\S 10]{P1}, then you might mull over why the previous example does not provide a counterexample.

\cmt{strictly speaking below, to apply Hamilton we would require curvature to be positive, but this is true instantly by the maximum principle here.}

\begin{figure}[h]
\begin{tikzpicture}

\newcommand\radius{0.5}

\draw [thick] (-5,\radius) -- (5,\radius);
\draw [thick] (-5,-\radius) -- (5,-\radius);


\draw (-5,\radius) arc[start angle=90,end angle=270, x radius=0.3*\radius, y radius=\radius];

\draw [dashed] (-5,-\radius) arc[start angle=-90,end angle=90, x radius=0.3*\radius, y radius=\radius];

\draw (5,\radius) arc[start angle=90,end angle=270, x radius=0.3*\radius, y radius=\radius];

\draw [dashed] (5,-\radius) arc[start angle=-90,end angle=90, x radius=0.3*\radius, y radius=\radius];


\draw [thick] (-5,\radius) arc[start angle=90,end angle=270, x radius=\radius, y radius=\radius];

\draw [thick] (5,-\radius) arc[start angle=-90,end angle=90, x radius=\radius, y radius=\radius];


\draw  [<->] (-5,-2*\radius) -- (0,-2*\radius) node[below]{$2$} -- (5,-2*\radius);
\draw  (-5,-2*\radius-0.1) -- (-5,-2*\radius+0.1);
\draw  (5,-2*\radius-0.1) -- (5,-2*\radius+0.1);



\draw [fill=black] (-0.1,\radius) arc[start angle=-180,end angle=0, x radius=0.1, y radius=0.08] node[above]{$p$};

%
%
%
%
%
%

\end{tikzpicture}
\caption{Thin cylinder with caps}
\label{fig1}
\end{figure}
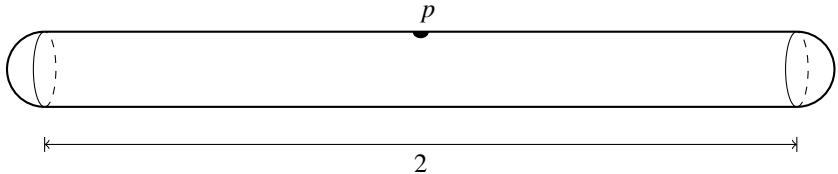

To construct Example \ref{ex3}, consider a 2-sphere that geometrically looks like a cylinder of length $2$ and small circumference, with small hemispherical caps at the ends. 
See Figure \ref{fig1}. We make the circumference small so that the total area is  $8\pi\ep\ll 1$. Thus the circumference of the cylinder is a little below $4\pi\ep$.
A theorem of Hamilton \cite{Ham88}  tells us that the subsequent Ricci flow $\ti g(t)$ will exist for a time $\ep$, and at that time the volume converges to zero and the curvature blows up everywhere. Just before time $\ep$, the surface looks like a tiny round sphere.
Certainly it is easy to see that the flow cannot exist any longer than this time since we can compute that
$$\frac{d}{dt}\Vol_{\tilde g(t)}S^2=\int_{S^2}\half\tr\pl{\ti g}{t} dV_{\ti g}
=-2\int_{S^2} K_{\ti g}\,dV_{\ti g}=-8\pi,$$
by Gauss-Bonnet.
To obtain the Ricci flow $g(t)$ on $D$ of the example, we pull back 
$\ti g(t)$ under the exponential map $\exp_p$, defined with respect to the initial metric (not the time $t$ metric), restricted to the unit disc $D$ in $T_pM$, where $p$ is any point midway along the cylinder.

\begin{figure}[h]
\begin{center}
\def\svgwidth{300pt}
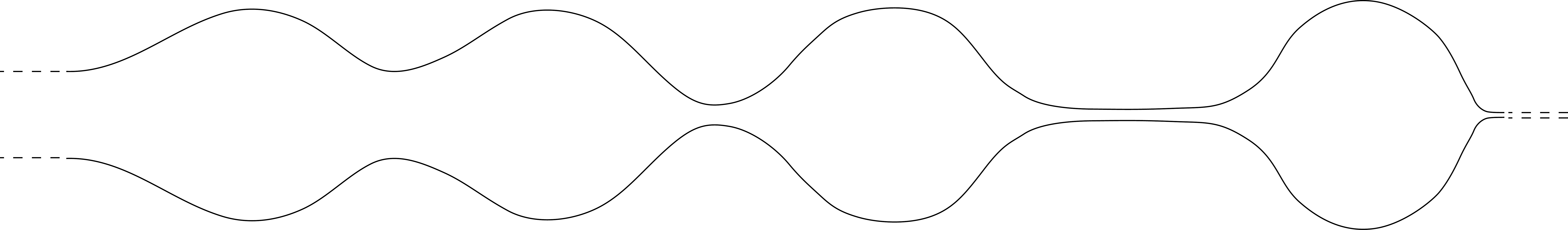
\end{center}
\caption{$S^2\times \R$ with thinner and thinner necks}
\label{bubble_fig}
\end{figure}

\begin{example}
\label{ex4}
Conjecturally, we have an example of a smooth complete three-manifold for which there  does not exist a smooth Ricci flow solution even for a short time. We construct this by connecting countably many three-spheres together by necks that progressively become longer and thinner. Each neck should want to pinch in a time that converges to zero as the necks get thinner and longer. 
See Figure \ref{bubble_fig} for a construction on $S^2\times \R$. Closing up the left-hand side would give an example on $\R^3$.
This construction can be made even under the assumption that the Ricci curvature is bounded below.
Making this example rigorous is currently beyond existing technology. One would need (for example) a pseudolocality theorem that is valid in the presence of unbounded curvature, but this is unavailable currently.
\end{example}

\subsection{What is known in the noncompact case}

There are several situations in which we can say a lot about Ricci flows even though the curvature can be unbounded. 


The first theory to be developed was the theory for two-dimensional underlying manifolds, which is also the one complex-dimensional K\"ahler case. This is the easiest case, and  gives the strongest results in terms of existence, uniqueness and asymptotic behaviour. The approach developed for this case has turned out to be important to understand the higher dimensional cases. We will see the main results in Section \ref{2Dsect}.

One way that higher dimensional Ricci flow can be made to have a clean existence theory for which there are no `unnecessary' singularities, just as there are no unnecessary singularities in the general two-dimensional case, is to impose a sufficiently strong positivity of curvature condition. In particular, Cabezas-Rivas and Wilking \cite{CR_W} developed such a theory in the case that the initial manifold is complete with nonnegative complex sectional curvature. 


Recently, several higher-dimensional cases that are substantially more singular have been treated. We will be looking particularly at the case in which the initial manifold is complete and three-dimensional, with a lower Ricci curvature bound.  Note that Example \ref{ex4} fits into this category, so we need to find a novel way of flowing. The solution is to flow locally, echoing the (quite different) work of Deane Yang \cite{DY} and the concepts of Ricci flow on manifolds with boundary (see \cite{panagiotis} and the earlier works referenced therein).
Once a solution has been constructed, and we have strong-enough control on the evolution of the Ricci curvature and the distance function, 
we will obtain applications to the theory of Ricci limit spaces.
See our work with Miles Simon \cite{ST2}, \cite{hochard} and Section \ref{localRFsect} for the results, and Section \ref{ACCTsect} for some applications.

In \cite{B_CR_W}, Bamler, Cabezas-Rivas and Wilking describe a  heat kernel method for controlling the curvature in certain cases, which allows the theory of \cite{CR_W} mentioned above to be generalised from nonnegative curvature to curvature bounded below, provided one assumes a global noncollapsing hypothesis. In \cite{lai}, this work was combined with that of \cite{ST2} to generalise both works. See also \cite{MT2}.

Finally, another interesting case in which to consider unbounded curvature is that of K\"ahler Ricci flow, though we will not attempt to survey that field here other than the one complex-dimensional case.


\subsection{The two-dimensional theory}
\label{2Dsect}

In 2D we will be able to flow any smooth initial metric. We do not require the curvature to be bounded, so Shi's theorem \ref{shi_thm} does not help. We do not even require that the initial metric is complete. 
Even if we can settle the existence problem, the level of generality here poses a potentially serious  problem of nonuniqueness of solutions, and to understand this let us consider the PDE more carefully.

\cmt{Saw the 2D case very early in the lectures - before Hamilton's Theorem}

We already saw that in 2D, the Ricci flow is the equation 
$$\pl{g}{t}=-2Kg,$$
where $K$ is the Gauss curvature of $g$. If we choose local isothermal coordinates $x,y$, i.e. so the metric $g$ can be written 
$e^{2u}(dx^2+dy^2)$ for a locally-defined scalar function $u$, then we can write the Gauss curvature as $K=-e^{-2u}\lap u$, where 
$\lap=\pll{}{x}+\pll{}{y}$ is only locally defined. 
Thus the Ricci flow can be written locally as
\beql{conf_factor_flow}
\pl{u}{t}=e^{-2u}\lap u,
\eeq
which is the so-called logarithmic fast diffusion equation (up to a change of variables). 
Given such a parabolic equation, the accepted wisdom is that one should specify initial and boundary data in order to obtain existence and uniqueness. However here we abandon the idea of specifying boundary data in the traditional sense and instead replace it with the the global condition that the flow should be complete for all positive times. 
Miraculously this completeness condition is just weak enough to permit existence, and just strong enough to force uniqueness.

\cmt{add ref to UK-Japan winter school?}

\begin{theorem}
\label{2D_eu}
Given any smooth (connected) Riemannian surface $(M,g_0)$, 
possibly incomplete and possibly with unbounded curvature,
define $T\in (0,\infty]$ depending on the conformal type of $(M,g_0)$ by
\[ T := \begin{cases}
 \frac1{4\pi}\Vol_{g_0}M & \text{if }(M,g_0)\cong \mathbb C\text{ or }\mathbb R P^2, \\
 \frac1{8\pi}\Vol_{g_0}M & \text{if }(M,g_0)\cong\mathcal S^2, \\
 \qquad\infty & \text{otherwise}.
\end{cases} \]
Then there exists a Ricci flow $g(t)$ on $M$ for $t\in [0,T)$ such that
\begin{enumerate}
\item
$g(0)=g_0$, and 
\item
$g(t)$ is complete for all $t\in (0,T)$.
\end{enumerate}
If $T<\infty$ then $\Vol_{g(t)}M\to 0 $ as $t\upto T$.
Moreover, this flow is unique in the sense that if $\ti g(t)$ is any other complete Ricci flow on $M$ with $\ti g(0)=g_0$, existing 
for $t\in [0,\ti T)$, then $\ti T\leq T$
and $\ti g(t)=g(t)$ for all $t\in [0,\ti T)$.
\end{theorem}
%
%
The existence of this theorem was proved with Giesen in \cite{GT2}, following earlier work in \cite{JEMS}. The compact case, where the metric is automatically complete and the curvature bounded, was already proved by Hamilton and Chow \cite{Ham88, chow}. Other prior results include the extensive theory of the logarithmic fast diffusion equation, typically posed on $\R^2$, for example \cite{DdP95, DD, vaz}. The uniqueness assertion was proved in 
\cite{ICRF_uniq} following a large number of prior partial results (see \cite{ICRF_uniq} for details).

We illustrate the theorem with two simple examples.

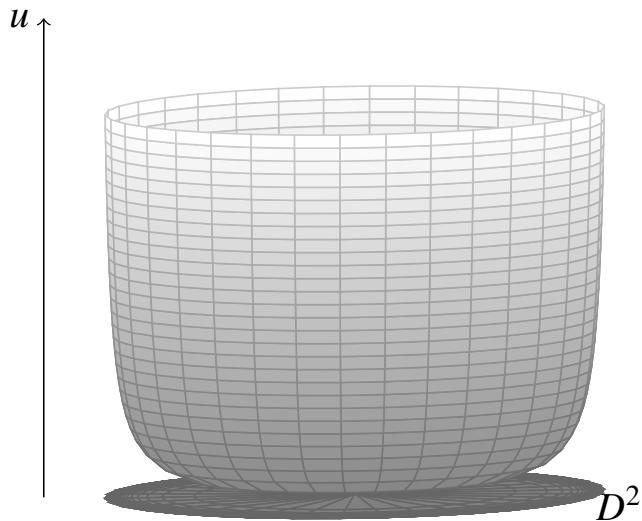
\begin{figure}[h]
\begin{tikzpicture}[scale=1.5]
\begin{axis}[hide axis, view/v=5,
  domain=1:20, 
  y domain=0:2*pi, 
  colormap={CM}{color=(gray) color=(white)}, 
  samples=30,
  samples y=35,
  z buffer=sort,
 ]
\addplot3[surf] 
  	({(1-1/x)^(1/2)*cos(deg(y))},{(1-1/x)^(1/2)*sin(deg(y))},{0.8/2})
	node[right]{\ $D^2$};  
\addplot3[surf] 
  	({(1-1/x)^(1/2)*cos(deg(y))},{(1-1/x)^(1/2)*sin(deg(y))},{x/2});
\addplot3[draw, ->] coordinates{(-1.4,0,0) (-1.4,0,12)} node[left] {$u$};
\end{axis}

\end{tikzpicture}
\caption{Metric for positive time stretches at infinity in Example \ref{discex}}
\label{disc_fig}
\end{figure}

\begin{example}
\label{discex}
Suppose $(M,g_0)$ is the flat unit disc. Finding a Ricci flow starting with $(M,g_0)$ is equivalent to solving \eqref{conf_factor_flow} with the initial condition $u|_{t=0}=0$. 
There is the obvious solution that has $u\equiv 0$ for all time, and we also saw the solution of Example \ref{ex3}. In general, standard PDE theory tells us that if we specify  boundary data on the disc as a function of time (and space) then there exists a unique solution that approaches that data at the boundary. 
None of these solutions give  complete metrics, so Theorem \ref{2D_eu} is giving another solution that blows up at the boundary. For this instantaneously complete flow, at small positive time $t>0$ the flat metric gets adjusted asymptotically near the boundary to the Poincar\'e metric scaled to have constant curvature $-\frac{1}{2t}$, but remains almost flat on the interior.
See Figure \ref{disc_fig}.
\end{example}



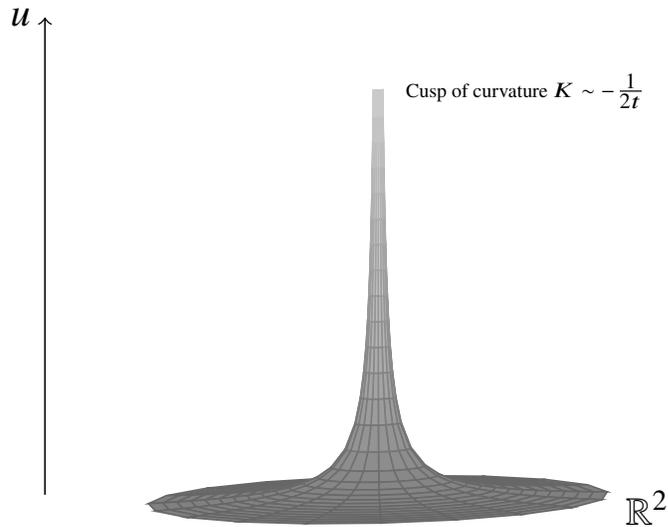
\begin{figure}[h]
\begin{tikzpicture} [scale=1.5]

\begin{axis}[hide axis, view/v=5,
  domain=1:7, 
  y domain=0:2*pi, 
  colormap={CM}{color=(gray) color=(white)}, 
  samples=15, 
  samples y=20, 
  z buffer=sort,
 ]

\addplot3[surf] 
  	({x*cos(deg(y))},{x*sin(deg(y))},{1/x})
	node[right]{\quad $\R^2$};  	
\addplot3[surf] 
  	({(1/x)*cos(deg(y))},{(1/x)*sin(deg(y))},{x})
	node[right]{\ \tiny Cusp of curvature $K\sim -\frac{1}{2t}$};  
\addplot3[draw, ->] coordinates{(-12,0,0) (-12,0,8)} node[left] {$u$};
\end{axis}

\end{tikzpicture}
\caption{Puncture turns into a hyperbolic cusp in Example \ref{puncex}}
\label{punc_fig}
\end{figure}

\begin{example}
\label{puncex}
Alternatively, suppose that $(M,g_0)$ is the flat Euclidean plane with one point removed to make it incomplete. This time the instantaneously complete solution stretches the metric immediately near the puncture to make it asymptotically like a hyperbolic cusp scaled to have curvature 
$-\frac{1}{2t}$.
See Figure \ref{punc_fig}.
\end{example}
%

It is important to appreciate that the flows themselves in Theorem \ref{2D_eu} can have unbounded curvature at spatial infinity at each time \cite{GT3}
as well as initially. 
More precisely, the curvature can be unbounded above; it follows from \cite{strong_uniqueness} that our solutions always satisfy the Gauss curvature estimate $K_{g(t)}\geq -\frac{1}{2t}$.
As alluded to in Remark \ref{Shi_length}, they might alternatively start as complete metrics of bounded curvature (so Shi's theorem \ref{shi_thm} applies), develop blow-up of curvature in finite time (at spatial infinity) but then continue smoothly for all time \cite{GT4}.
These types of behaviour can even happen under strong positivity of curvature conditions, if one permits higher dimensional flows
\cite{CR_W}.

\cmt{missing out maximally stretched, asymptotic behaviour, 
Schwarz-Pick-Ahlfors etc}

\subsection{Proof ideas for the 2D theory}

It will be helpful to outline some of the ideas in the proof of the existence of Theorem \ref{2D_eu}; we'll emphasise the ones that 
ended up being useful in the higher dimensional theory.
A rough principle is that we would like to bring Shi's result into play despite the unboundedness of the curvature. We achieve this through the following steps, where we are free to assume that $M$ is not closed (or we would be in the classical situation of Hamilton).

\begin{enumerate}
\item
First we take an exhaustion of $M$ by sets $\Om_1\subset\Om_2\subset \cdots\subset M$, each compactly contained in $M$ and with smooth boundary.
\item
Next, for each $i\in\N$ we restrict $g_0$ to $\Om_i$ and blow up the metric conformally near the boundary to make the metric complete and of very large constant negative curvature near the boundary. Call the resulting metric $g_i$.
\item
Since each $(\Om_i,g_i)$ is complete, with bounded curvature, we can apply Shi's theorem \ref{shi_thm} to obtain Ricci flows $g_i(t)$.
\item
Next, we prove that each $g_i(t)$ exists for a uniform time (in fact for 
\emph{all} time in this case) and enjoys uniform curvature estimates (independent of $i$).
\item
The uniform estimates from the previous step give us compactness: a subsequence converges to a limit flow $g(t)$, this time on the whole of $M$.
\item
Finally we argue that the limit flow $g(t)$ is complete for all positive times, and has $g_0$ as initial data.
\end{enumerate}
We'll reuse the same strategy for higher dimensions, but the results are a little different.

%
%

\section{Ricci curvature and Ricci limit spaces}


Most of the remaining material we will cover in these lectures will be geared to understanding Ricci flow in higher dimensions in the unbounded curvature case. The new ideas we  wish to explain are best illustrated in the case of three-dimensional manifolds that have a lower Ricci curvature bound. These ideas generalise to higher dimensions if we assume correspondingly stronger curvature bounds. This theory has a number of potential applications to the understanding of the topology and geometry of manifolds with lower curvature bounds, but the application on which we focus in these lectures will be to so-called Ricci limit spaces, which we shall shortly define as some sort of limits of Riemannian manifolds with lower Ricci bounds. 
Before that, we survey some of the most basic consequences of a manifold having a  lower Ricci bound.

\subsection{Volume comparison}
\label{comp_geom_sect}

\cmt{We will need at some point to be able to say that we can go from 
\eqref{conditions2} to the hypotheses of the extension lemma.
New $v_0$ depends on $r_1$ in addition to $\al_0$ and $v_0$.
But MOST ATTENTION required to fact that it depends on $r$...}

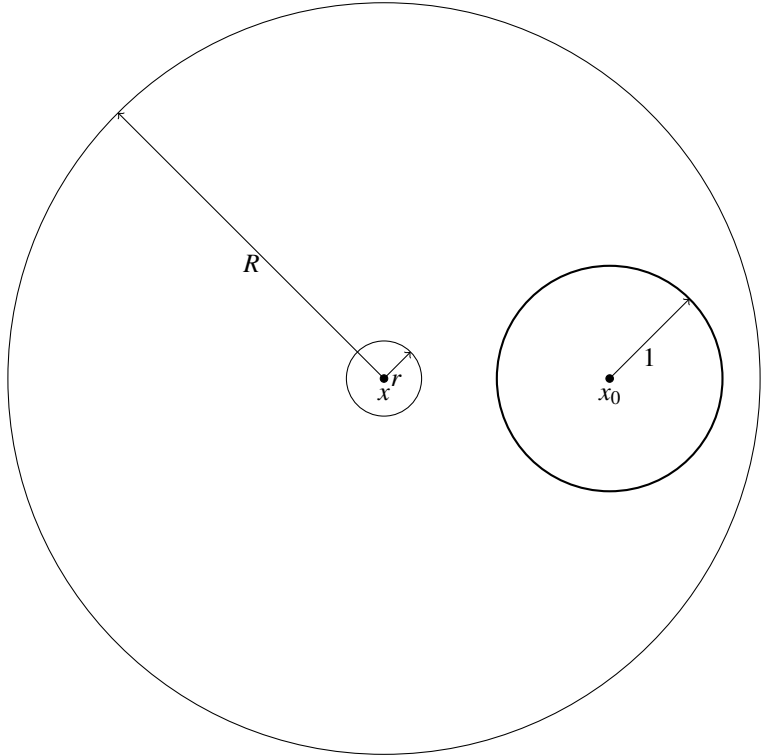
\begin{figure}[h]
\centering
\begin{tikzpicture}[xscale=0.5,yscale=0.5]

\newcommand\sq{0.7071068} 

\draw (0,0) circle (10cm);
\draw (0,0) circle (1cm);
\draw [thick] (6,0) circle (3cm);
\draw [fill] (6,0) circle (0.1cm) node[below]{$x_0$};
\draw [fill] (0,0) circle (0.1cm) node[below]{$x$};

\draw [->] (0,0) -- (\sq,\sq) node[midway, below]{$r$};
\draw [->] (6,0) -- (6+3*\sq,3*\sq) node[midway, below]{$1$};
\draw [->] (0,0) -- (-10*\sq,10*\sq) node[midway, below]{$R$};

\end{tikzpicture}
\caption{Comparing the volumes of nonconcentric balls.}
\label{comparison_balls_fig}
\end{figure}

Ricci curvature controls the growth of the volume of balls as the radius increases. The simplest result of this form is:

\begin{theorem}[Bishop-Gromov, special case.]
Given a complete Riemannian $n$-manifold $(M,g)$ with $\Ric\geq 0$, and a point $x_0\in M$, the \emph{volume ratio} function
\beql{vol_ratio}
r\mapsto \frac{\VolB_g(x_0,r)}{r^n}
\eeq
is a weakly decreasing (i.e. nonincreasing) function of $r>0$.
\label{BGthm}
\end{theorem}
Since the limit of the function \eqref{vol_ratio} as $r\downto 0$ is necessarily equal to $\om_n$, the volume of the unit ball in Euclidean $n$-space, this implies that 
$\VolB_g(x_0,r)\leq \om_n r^n$. 
But of course, the theorem also gives us \emph{lower} volume bounds:
For $0<r\leq R$, we have
$$\VolB_g(x_0,r)\geq (r/R)^n \VolB_g(x_0,R).$$
In fact, we do not even need to have the same centres here. Later it will be useful to observe that for any $x\in M$ we can set
$R=r+1+d(x,x_0)$ and estimate,
as in Figure \ref{comparison_balls_fig}, that
\beql{vol_rat_lower_bd}
\frac{\VolB_g(x,r)}{r^n}\geq R^{-n} \VolB_g(x,R)
\geq R^{-n} \VolB_g(x_0,1),
\eeq
and thus knowing that one unit ball $\VolB_g(x_0,1)$ has a specific lower volume bound then gives a lower bound for the volume of any other ball, though that lower bound gets weaker as $x$ drifts far from $x_0$.

\cmt{but shouldn't be sharp. Once $\VolB_g(x_0,1)$ is fixed (smallish) and $R$ is fixed (large) then the smallest we can make $\VolB_g(x,r)$ is at the tip of a cone... Suggests lower bound like $R^{1-n}$
in this $\Ric\geq 0$ case.}

\cmt{In fact, we'll need the analogous argument when the Ricci lower bound is negative, not just zero.
But the easiest way of explaining this seems to be to give the zero lower bound case, as here, and then later to say that it generalises OK. Note that we need this when doing the extension lemma, but also when doing Gromov compactness}

Theorem \ref{BGthm} also allows us to define a finite \emph{asymptotic volume ratio}
$${\rm AVR}:=\lim_{r\to\infty}\frac{\VolB_g(x_0,r)}{r^n}.$$ 
In practice it is often particularly relevant whether AVR equals zero or is positive.

\brmk
\label{gen_al0_rmk}
The Bishop-Gromov theorem as stated extends easily to manifolds with more general Ricci lower bounds $\Ric\geq -\al_0$, but now the volume of a ball of radius $r$
can grow as fast as it does in the complete model space  of constant sectional curvature where the constant is chosen so that 
the Ricci curvature equals $-\al_0$.
\ermk

\cmt{A more general form of the comparison compares the volumes of annuli in the space and the model space.
Let $0\leq r_1<r_2,r_3\leq r_4$. Then 
$$\frac{\Vol(A_{r_3,r_4}(p))}{\Vol(A_{r_3,r_4}^H)}\leq
\frac{\Vol(A_{r_1,r_2}(p))}{\Vol(A_{r_1,r_2}^H)}.$$
}



\subsection{What is a Ricci limit space?}
\label{what_is_RLS_Sect}

Suppose that $(M_i^n,g_i,x_i)$ is a sequence of pointed complete Riemannian manifolds, i.e. a sequence of complete Riemannian $n$-manifolds 
$(M_i^n,g_i)$ together with points $x_i\in M_i$,
and suppose that the sequence satisfies the properties that
\beql{conditions}
\left\{
\begin{aligned}
\Ric_{g_i}\geq -\al_0\\
\VolB_{g_i}(x_i,1)\geq v_0>0.
\end{aligned}
\right.
\eeq
%
%
Such a sequence of manifolds with lower Ricci curvature control has some sort of compactness property. That is, a subsequence will converge in a suitable sense. What it converges to will be more general than a Riemannian manifold, but will be a very special type of metric space. These limits will be what we mean by Ricci limit space. Before giving any details, let's imagine an example.

\begin{figure}[h]
\centering
\begin{tikzpicture}

\newcommand\xtrans{0}
\draw [thick] (\xtrans+-1,-2) -- (\xtrans+0,0) -- (\xtrans+1,-2);
\draw [dashed] (\xtrans+0.5,-1) arc[start angle=0,end angle=180, 
x radius=0.5, y radius=0.25];
\draw [thick] (\xtrans-0.5,-1) arc[start angle=180,end angle=360, 
x radius=0.5, y radius=0.25];

\newcommand\smsc{0.5} 
\renewcommand\xtrans{-7}
\draw [thick] (\xtrans+-1,-2) -- (\xtrans-\smsc/2,-\smsc);
\draw [thick] (\xtrans+\smsc/2,-\smsc) -- (\xtrans+1,-2);
\draw[thick, domain=(\xtrans-\smsc/2):(\xtrans+\smsc/2), samples=50] 
plot (\x, {
-2*(\x-\xtrans)*(\x-\xtrans)/\smsc-\smsc/2
});
\draw [dashed] (\xtrans+0.5,-1) arc[start angle=0,end angle=180, 
x radius=0.5, y radius=0.25];
\draw [thick] (\xtrans-0.5,-1) arc[start angle=180,end angle=360, 
x radius=0.5, y radius=0.25];

\renewcommand\smsc{0.2} 
\renewcommand\xtrans{-4}
\draw [thick] (\xtrans+-1,-2) -- (\xtrans-\smsc/2,-\smsc);
\draw [thick] (\xtrans+\smsc/2,-\smsc) -- (\xtrans+1,-2);
\draw[thick, domain=(\xtrans-\smsc/2):(\xtrans+\smsc/2), samples=50] 
plot (\x, {
-2*(\x-\xtrans)*(\x-\xtrans)/\smsc-\smsc/2
});
\draw [dashed] (\xtrans+0.5,-1) arc[start angle=0,end angle=180, 
x radius=0.5, y radius=0.25];
\draw [thick] (\xtrans-0.5,-1) arc[start angle=180,end angle=360, 
x radius=0.5, y radius=0.25];

\draw [thick, ->] (-2-0.2,-1) -- (-2+0.2,-1);

\end{tikzpicture}
\caption{Smoothed out cones approaching a Ricci limit space.}
\label{fig_cones}
\end{figure}
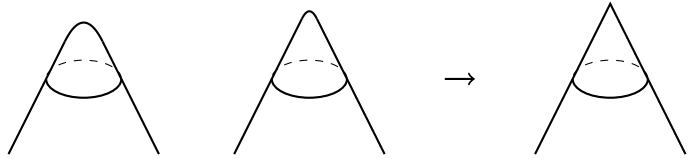

Consider the two-dimensional cone 
$\cc:=(\R^2\setminus\{0\},dr^2+\al r^2 d\th^2)$, where $\al\in (0,1)$ 
and we are working in polar coordinates.
Intuitively there is a delta function of positive curvature at the vertex, i.e. at the origin, scaled by a factor $2\pi(1-\al)$.
If we view $\cc$ as a metric space, and add in the origin, we have an example of a Ricci limit space with the origin as a singular point (these terms will be defined in a moment).
This space will be the limit of a sequence of surfaces $(M_i,g_i)$
that are constructed by smoothing out the cone point of $\cc$ by a smaller and smaller amount as $i$ gets large. Each of these surfaces can be assumed to have nonnegative Gauss curvature, and thus Ricci curvature bounded below. Moreover, if we let $x_i$ be the origin for each $i$, we can arrange that $\VolB_{g_i}(x_i,1)$ converges to the area of the unit ball in $\cc$ centred at the origin, i.e. to $\pi\al$, which is positive.
Thus the conditions \eqref{conditions} will be satisfied, and we see
that the cone $\cc$ is a Ricci limit space; see Figure \ref{fig_cones}.

So far this is very vague. The notion of convergence we are considering here is that of \emph{pointed Gromov-Hausdorff} convergence. Before defining that, it makes sense to recall
the usual notion of Gromov-Hausdorff convergence.
A sequence of metric spaces $(X_i,d_i)$ converges to a limit metric space $(X,d)$ in the 
Gromov-Hausdorff sense if  for all
$\ep>0$ and for sufficiently large $i$ (depending also on $\ep$) 
there exists a map $f:X_i\to X$ such that
\begin{enumerate}
\item
for all $x,y\in X_i$ we have
$$|d(f(x),f(y))-d_i(x,y)|<\ep,$$
\item
for every $z\in X$, there exists $x\in X_i$ such that
$$d(f(x),z)<\ep.$$
\end{enumerate}
This is clearly a very weak notion of closeness. For example, a finer and finer lattice of points in the plane would converge to the plane itself.

\cmt{not saying anything about uniqueness etc}

We'll loosely call any such map $f$ an $\ep$-Gromov-Hausdorff approximation. 

\cmt{nb. need pointed GH convergence in setting of sequences of metric spaces rather than manifolds when considering tangent cones. At least at first glance we do. Discussion there was adjusted to cover this issue.}

The general definition of \emph{pointed} Gromov-Hausdorff convergence is similar:
Given a complete metric space $(X,d)$ and a point $x_\infty\in X$, 
i.e. a pointed metric space $(X,d,x_\infty)$, 
we say that a sequence of pointed metric spaces
$(X_i,d_i,x_i)$ converges to $(X,d,x_\infty)$ in 
the pointed Gromov-Hausdorff sense if
given any radius $r>0$ and any $\ep>0$, then for sufficiently large $i$ (depending also on $r$ and $\ep$)
there exists a map $f:B_{d_i}(x_i,r)\to X$ such that
\begin{enumerate}
\item
$f(x_i)=x_\infty$,
\item
for all $x,y\in B_{d_i}(x_i,r)$ we have
$$|d(f(x),f(y))-d_i(x,y)|<\ep,$$
\item
the $\ep$-neighbourhood of the image $f(B_{d_i}(x_i,r))$
contains the ball $B_d(x_\infty, r-\ep)\subset X$.
\end{enumerate}
\cmt{pedantry corner - the condition as above is not quite as in BBI because the distortion is the SUP of what we estimate, so could EQUAL $\ep$}
In practice, we will essentially only consider this definition in the special case that the metric spaces $(X_i,d_i,x_i)$ are in fact Riemannian manifolds 
$(M_i^n,g_i,x_i)$. This restriction removes some 
annoying pathology by virtue of restricting us to working with length spaces and boundedly compact spaces (see \cite{BBI} for definitions).

\cmt{Do we get example of nonuniqueness if we have completely general pointed GH convergence?}

Indeed, working only with boundedly compact spaces, our limit $(X,d,x_\infty)$ is unique in the sense that any other (complete) limit will be isometric to this via an isometry that fixes the marked points 
\cite[Theorem 8.1.7]{BBI}.


Another property that follows with this restriction
to Riemannian manifolds, by virtue of restricting us to length spaces,   is that for each $r>0$, the balls $B_{g_i}(x_i,r)$ converge to the ball $B_d(x_\infty,r)$ in the Gromov-Hausdorff sense, which is not true in general.
Beware that the converse fails: Gromov-Hausdorff convergence of open balls of arbitrary radius $r$ does not imply the pointed convergence as given here. Nor does convergence of closed balls of arbitrary radius.
(Just consider metric spaces each consisting of two points only.)

\cmt{Example of not having GH converge of balls in general case: take a sequence of metric spaces each consisting of two point, separated by a distance $1-1/i$. Open unit balls do not converge.}

\cmt{Example of converse failing: take points in the plane at $(\pm 10, 1)$ and 
$(\pm 10, -2)$. Set $p=(-10,0)$ and $q=(10,0)$. 
Fix the induced distances, and then forget about the ambient space.
Now reduce a little the distance across one diagonal, say
$(-10,-2)$ to $(10,1)$. 
Balls centred at $p$ are isometric to balls centred at $q$, but not via an isometry that sends $p$ to $q$. This can be used as a building block. 
This does not cover the manifold version of this question.}

See \cite{BBI}, p.271\emph{ff} for further details on this general topic.

Let's return to the original sequence of pointed manifolds 
satisfying conditions \eqref{conditions}.
By virtue of the lower Ricci curvature bound, and the resulting volume comparison discussed 
in Section \ref{comp_geom_sect}, a compactness result of Gromov given in Theorem \ref{grom_cpt_thm} of Appendix \ref{grom_appen} tells us that we can pass to a  subsequence in $i$ and obtain 
pointed Gromov-Hausdorff convergence 
\beq
\label{GH_cgnce_def}
(M_i,g_i,x_i)\to (X,d,x_\infty)
\eeq
to some complete metric space 
$(X,d)$ (even a length space) and $x_\infty\in X$.

\begin{definition}
A complete metric space $(X,d)$ arising in a limit \eqref{GH_cgnce_def} of a sequence
$(M_i,g_i,x_i)$ of pointed complete Riemannian $n$-manifolds satisfying \eqref{conditions}
is called a (noncollapsed) Ricci limit space.
\end{definition}
In these lectures we will not consider so-called collapsed Ricci limit spaces, i.e. we will always assume the volume lower bound of \eqref{conditions}.

Such Ricci limit spaces were studied extensively starting with Cheeger and Colding around the late 90s. 
They can be considered as rough spaces that can be viewed as having some sort of lower Ricci bound. As such they are analogues of Alexandrov spaces, which are metric spaces with a notion of lower \emph{sectional} curvature bound. A synthetic notion of rough space with lower Ricci curvature bound (i.e. defined directly rather than as  a limit of smooth spaces) is provided by the so-called RCD spaces, for example; see \cite{villani_survey} for background.
The exact link between Ricci limit spaces and RCD spaces is a topic of intense current study. One of the basic questions in that direction will be 
settled by the Ricci flow methods we will see in the next section; see Remark \ref{RCDrmk}.


\subsection{How irregular can a Ricci limit space be?}
\label{irreg_sect}


Let's return to the example of a Ricci limit space that we have already seen, i.e. the two-dimensional cone $\cc$. Intuitively it has precisely one singular point; we now justify that by giving a definition of \emph{singular point}.
At each point $p$ in a Ricci limit space $(X,d)$, we can define a \emph{tangent cone} to be a (pointed Gromov-Hausdorff) 
limit of rescalings
$(X,\la_i d,p)$ for some $\la_i\to\infty$. We can always find a tangent cone by picking an arbitrary sequence $\la_i\to\infty$ and passing to a subsequence using Gromov compactness
as discussed in Theorem \ref{grom_cpt_thm}. 
Strictly speaking we have only talked about Gromov compactness for sequences of Riemannian manifolds rather than metric spaces as here, but $(X,\la_i d,p)$ itself is the pointed Gromov-Hausdorff limit of $(M_j,\la_i^2 g_j,y_j)$ as $j\to\infty$, for some sequence of points $y_j\in M_j$, and the tangent cone can be viewed as a limit of rescaled and rebased pointed manifolds
$(M_i,\mu_ig_i,y_i)$.

There may be different tangent cones corresponding to different sequences (see e.g. \cite[Theorem 8.41]{CC1}). 
Tangent cones have a lot of structure, see \cite{cheeger_book}
for details.

\begin{definition}
A point $p$ in a Ricci limit space $(X,d)$ is singular if at least one tangent cone is not Euclidean space.
\end{definition}
In fact, an elegant argument shows that if any one tangent cone is Euclidean space, then any other tangent cone (i.e. for any other sequence $\la_i\to\infty$) is Euclidean space too, see \cite{cheeger_book}.

In our cone example, the tangent cones are Euclidean space away from the vertex (thus regular) and a nontrivial cone at the vertex (hence singular).



\begin{figure}
\centering
\begin{tikzpicture}[xscale=2,yscale=2]

\newcommand\ang{50}


\draw [thick] ({sin(\ang)},{-cos(\ang)}) arc [radius = 1, start angle = -90+\ang, end angle = 270-\ang];

\draw [dashed] ({-sin(\ang)},{-cos(\ang)}) arc [radius = 1, start angle = 270-\ang, end angle = 270+\ang];

\draw [thick] ({sin(\ang)},{-cos(\ang)}) -- (0,{-1/(cos(\ang))});

\draw [thick] ({-sin(\ang)},{-cos(\ang)}) -- (0,{-1/(cos(\ang))});


\draw [dashed] ({sin(\ang)},{-cos(\ang)}) arc [x radius={sin(\ang)}, y radius = {0.1*(sin(\ang))}, start angle = 0, end angle = 180];

\draw ({-sin(\ang)},{-cos(\ang)}) arc [x radius={sin(\ang)}, y radius = {0.1*(sin(\ang))}, start angle = 180, end angle = 360];


\draw [dashed] (0,0) ellipse (1 and 0.1);

\end{tikzpicture}
\caption{Ice cream cone}
\label{fig_ice_cream_cone}
\end{figure}
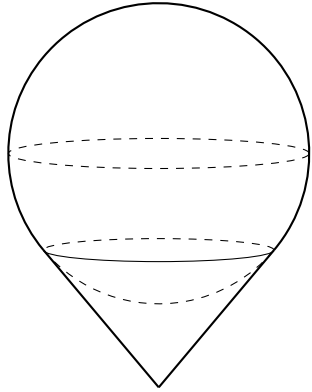

The singular set could be a lot worse in size, for example it could be dense. Let's see at least how one can construct a singular set with an accumulation point. To do this, start with a round sphere. Select a spherical cap that is strictly smaller than a hemisphere, and replace it with the top of a cone so that the surface is $C^1$, i.e. the tangent spaces match up. Thus we end up with an `ice cream cone' as in Figure \ref{fig_ice_cream_cone}.


We can now pick a countable collection of further spherical caps, of size converging to zero, that are mutually disjoint. Each of them can be replaced by a cone as before. As the caps become smaller, the cones become `flatter'. Each cone vertex will be a singular point. Clearly we have constructed another Ricci limit space, this time with accumulating singular points, which we call $(X_0,d_0)$.

Although the singular set here is infinite, it nevertheless has Hausdorff dimension zero, which is optimal:

\begin{theorem}[{Cheeger-Colding \cite{CC1}.}]
A noncollapsed Ricci limit space has Hausdorff dimension $n$ and its singular set 
has Hausdorff dimension no larger than $n-2$.
\label{CCthm1}
\end{theorem}
\cmt{Aside: with 2-sided bounds, these cones are ruled out, and the dimension of the singular set is codimension 4.}
By considering the product of $(X_0,d_0)$ with $\R^{n-2}$, which is clearly another Ricci limit space, we see that the control on the size of the singular set of Theorem \ref{CCthm1} is optimal in all dimensions.

Our example $(X_0,d_0)$ may have an infinite singular set, but it is \emph{topologically} regular in the sense that it is \emph{homeomorphic} to a smooth manifold, in this case $S^2$. The simple conical points present in this example have neighbourhoods homeomorphic to flat discs.
However even this milder form of regularity fails if we increase the dimension from $2$ to $4$ because of the following example.

\begin{example}
\label{eguchi_hanson}
There exists a four-dimensional complete Riemannian manifold $(M,g_0)$, known as the Eguchi-Hanson manifold, 
that is Ricci-flat, and which looks asymptotically like $\R^4/\sim$,
where $\sim$ is the equivalence relation that identifies each point $x\in\R^4$ with $-x$. By picking an arbitrary point $x_0\in M$ and 
taking pointed manifolds $(M_i,g_i,x_i):=(M,\frac{1}{i}g_0,x_0)$  for each $i$, which satisfy the conditions \eqref{conditions}, and converge in the pointed Gromov-Hausdorff sense to $\R^4/\sim$, we find that $\R^4/\sim$ is a Ricci limit space. But this is a cone over $\R P^3$, and being a cone over 
something other than a sphere, it cannot be homeomorphic to a manifold. 
\end{example}

\brmk
\label{cone_remark}
The term \emph{cone} can refer to many different things.
Here when we say cone over $\R P^3$ we can start with the Riemannian cone defined on $(0,\infty)\times\R P^3$, with the warped product metric
$dr^2+r^2 g_1$, where $r\in (0,\infty)$ and $g_1$ is the standard metric on $\R P^3$,
then view the space as a metric space and complete it by adding one point at the tip of the cone.
\ermk
\cmt{Amazingly we seem to have got away without defining a metric cone - only ever talking about Riemannian cones.}
\cmt{If homeomorphic, then it's still homeo after removing the vertex.
Thus we have that $S^3\times I$ homeo to a subset of link times $I$. 
This forces link to be simply connected. }
By taking a product of the Eguchi-Hanson manifold and $\R^{n-4}$, 
we see that in general the best regularity one can hope for is for our Ricci limit space to be homeomorphic to a manifold off a set of dimension $n-4$. 
It is a long-standing conjecture, normally attributed to Anderson, Cheeger, Colding and Tian (ACCT) that in some sense this is the case.
We'll give a precise formulation of the conjecture in general dimension in a moment, once we have introduced some more theory. 
But we can immediately state it precisely in the three-dimensional case: The whole 
Ricci limit space should be homeomorphic to a topological manifold.
With Miles Simon we prove more:

\begin{theorem}[{with Miles Simon \cite{ST2}.}]
In the case $n=3$, any Ricci limit space is locally bi-H\"older homeomorphic to a smooth manifold.
\label{ACCT3Dthm}
\end{theorem}
The proof of this will require a string of new ideas, some of which are
due to Hochard \cite{hochard}. In fact, a slight improvement is possible.

\begin{theorem}[{with Andrew McLeod \cite{MT}, using mainly the theory from 
\cite{ST1, ST2}.}]
In the case $n=3$, any Ricci limit space is globally homeomorphic to 
a smooth manifold via a homeomorphism that is 
locally bi-H\"older.
\label{pyramid_version}
\end{theorem}
Prior to these results, it was only possible to give topological regularity statements about the \emph{regular set} of a Ricci limit space (i.e. the complement of the singular set) or suitable generalisations thereof. This reduces considerations to spaces that are almost-Euclidean in some sense, which is a completely different task. However, it does have the advantage of having a chance of working in general dimension.

\begin{theorem}[{Cheeger-Colding \cite{CC1}.}]
There exists an open neighbourhood of the regular set of a general Ricci limit space that is locally bi-H\"older homeomorphic to a smooth manifold.
\label{CCbiholder}
\end{theorem}
In Theorem \ref{CCbiholder} we can make the H\"older exponent as close as we like to $1$ by shrinking the set (still containing the regular set). It is unknown whether we can make the homeomorphism Lipschitz.
\cmt{see p. 411 of \cite{CC1}.}

\cmt{give the definition of $\calr_\ep$ etc? Quite a bit of side explanation would be required to justify this. For example, the theory justifying that almost positive Ricci and almost-maximal volume propagates to smaller scales. And maybe volume convergence.}

We now return to formulate the ACCT conjecture in general dimension.
For each $k\in \{0,\ldots, n-1\}$, consider the subset $\s_k$ 
of the singular set $\s$ consisting of points at which no tangent cone can be written as a product of $\R^{k+1}$ with something else.
(At most they can be a product of $\R^k$ with something else.)
Note that 
$$\s_0\subset\s_1\subset\cdots\subset\s_{n-1}= \s.$$
Cheeger-Colding's theory \cite{CC1} tells us that the Hausdorff dimension of $\s_{k}$ is at most $k$. As usual, we are considering noncollapsed Ricci limit spaces \cite{CC1}. In fact, in this case, Cheeger-Colding also proved that $\s_{n-2}=\s$, which is how one proves Theorem \ref{CCthm1}.

The full ACCT conjecture is:

\begin{conjecture}[{\cite[Conjecture 0.7]{CC1}.}]
Given any Ricci limit space $X$, the interior of $X\setminus \s_{n-4}$
is homeomorphic to a topological manifold.
\end{conjecture}


\cmt{It is important to give the `more significantly' remark below because otherwise we could just take the doubling of the unit ball in the cone as the RCD space example.
(Or the related spherical suspension.)
Then it is clear that any global approximators must be closed Riemannian manifolds and the result follows for softer reasons - the space of smooth closed approximators is heavily restricted, and the technology of Ricci flow on closed manifolds can be invoked to give the restricted result}

\cmt{Got to be careful with the wording here because of the usual problems of making sense of the distance function on incomplete manifolds - looks OK though.
}

\brmk
\label{RCDrmk}
One consequence of Theorem \ref{ACCT3Dthm} is that we can settle the most obvious questions concerning whether synthetically-defined metric spaces with lower Ricci bounds 
always arise as noncollapsed Ricci limit spaces. 
One can see the cone over $\R P^2$ (plus the vertex) as such a metric space with non-negative Ricci curvature (e.g. a $\mathrm{RCD}^*(0,3)$ space, as shown in \cite{ketterer})
but it is not a manifold, being a cone over something other than a sphere, and thus cannot be a Ricci limit space by our theorem. More significantly, our theory in \cite{ST2} gives a local distinction in that a unit ball centred at the vertex of the cone cannot arise as the pointed Gromov-Hausdorff limit of a sequence of possibly incomplete Riemannian manifolds $(M^3_i,g_i,x_i)$ each with the property that $B_{g_i}(x_i,1)\subset\subset M_i$ (i.e. each manifold can fit in a unit ball) and $\Ric_{g_i}\geq 0$, say.
Note that it is much easier (and much weaker) to show that the cone over $\R P^2$ cannot arise as the pointed Gromov-Hausdorff limit of a sequence of \emph{closed} Riemannian 3-manifolds each with $\Ric_{g_i}\geq 0$. A global restriction of this form imposes very strong rigidity on the geometry and topology of the approximating closed manifolds. 
It is not even clear that a \emph{smooth} complete Riemannian 3-manifold with non-negative Ricci curvature can be approximated by closed 3-manifolds with the same curvature restriction.
\ermk

\subsection{How to apply Ricci flow to Ricci limit spaces}
\label{how_to_apply_sect}

The central question when trying to apply Ricci flow to Ricci limit spaces is: 
\begin{quote}
Can we start the Ricci flow with a Ricci limit space as initial data?
\end{quote}
%
Clearly we will have to revisit the traditional viewpoint of how to pose Ricci flow that we described at the start of the lectures. Now, when we are given initial data that is a metric space $(X,d)$ rather than a Riemannian manifold, we are asking Ricci flow to create not only a one-parameter family of metrics $g(t)$, but also the underlying manifold on which they live! The family $g(t)$ will now only live for positive times (i.e. not at $t=0$) and we will want $(M,d_{g(t)})$ to converge in a sufficiently strong sense to the initial data $(X,d)$ as $t\downto 0$.

Slightly more precisely, we will want that the distance metric $d_{g(t)}$ converges nicely to a distance metric $d_0$ on $M$ as $t\downto 0$, 
and that $(M,d_0)$ is isometric to the Ricci limit space $(X,d)$.
Here \emph{nice} convergence will be more than local uniform convergence on the domain $M\times M$ of the distance metrics.

Luckily we are not trying to flow a general metric space, but instead a Ricci limit space, which has a lot more structure. We can therefore see a potential strategy that essentially follows the strategy that 
worked in the 2D case. (We will need more ideas from the 2D case in due course.)
For a Ricci limit space $(X,d)$ arising as the pointed limit of 
$(M_i,g_i,x_i)$ we might try to:


\begin{enumerate}
\item
Flow each $(M_i,g_i)$ to give Ricci flows $(M_i,g_i(t))$, and argue that they exist on a time interval $[0,T]$ that is independent of $i$. ($T$ is allowed to depend on $\al_0$ and $v_0$.)
\item
Derive uniform estimates on the full curvature tensor of $g_i(t)$ that are independent of $i$, but can depend on $t\in (0,T]$ as well as 
$\al_0$ and $v_0$.
\item
Use parabolic regularity to improve these bounds to bounds on all derivatives of the full curvature tensor.
\item
Appeal to the smooth `Cheeger-Gromov-Hamilton' compactness that one obtains from these curvature bounds. More precisely, after passing to a subsequence, 
there exist a smooth manifold $M$, a smooth Ricci flow $g(t)$ for $t\in (0,T]$ and a point $x_\infty\in M$ such that
$$(M_i,g_i(t),x_i)\to (M,g(t),x_\infty)$$
smoothly as $i\to\infty$.
Smooth convergence means that
there exist an exhaustion $x_\infty\in\Om_1\subset\Om_2\subset\cdots\subset M$ 
and diffeomorphisms $\vph_i:\Om_i\to M_i$ such that $\vph_i(x_\infty)=x_i$ and 
$$\vph_i^*g_i(t)\to g(t)$$
smoothly locally on $M\times (0,T]$.
See \cite[\S 7]{RFnotes} for a more detailed description of Cheeger-Gromov convergence.
\item
Show that the distance metric $d_{g(t)}$ converges to some limit metric $d_0$ at least uniformly as $t\downto 0$.
(In practice, stronger convergence will be required.)
\item
Show that having passed to a subsequence, the sequence $(M_i,g_i,x_i)$ converges 
in the pointed Gromov-Hausdorff sense to $(M,d_0,x_\infty)$.
But this sequence also converges to the Ricci limit space we are trying to flow, so that must be isometric to $(M,d_0,x_\infty)$. 
\end{enumerate}
Note that the smooth manifold $M$ that we require has been 
\emph{created} by the Ricci flow.

If this programme could be completed, then we would have succeeded in starting the Ricci flow with a Ricci limit space.
Additionally, if we could prove that the distance metric $d_{g(t)}$ is sufficiently close to $d_0$, depending on $t$, in a strong enough sense, we can hope that the identity map from $(M,d_0)$ to $(M,d_{g(t)})$, i.e. to the smooth Riemannian manifold $(M,g(t))$, will be bi-H\"older. This would give a proof of Theorem \ref{ACCT3Dthm}.

\cmt{Not talking about bi-Holder either to original approximating metrics.}

Unfortunately, this strategy is a little naive as written. Even the first step fails at the first  hurdle: Example \ref{ex4} showed that we cannot expect to even start the Ricci flow with a general smooth complete manifold even if the Ricci curvature is bounded below.
Miles Simon carried the programme through in the case that the approximating manifolds are closed \cite{simon2012} so this problem doesn't arise.
For the general case this problem forces us to work locally.
Once we do that, the strategy will work.

\cmt{maybe say that In Section \ref{localRFsect} we'll consider the local Ricci flow theory that will allow us to carry out 1,2,3 above, then we'll state the RF from a RLS theorem in Section \ref{ACCTsect}.}

\section{Local Ricci flow in 3D}
\label{localRFsect}

\subsection{Local Ricci flow theorem with estimates}

The following theorem, proved with Miles Simon,  gives the local Ricci flow that we need. We can apply it to each of the 
$(M_i,g_i,x_i)$ considered in the previous section in order to make the strategy there work.
As mentioned earlier, local Ricci flow under different hypotheses or with weaker conclusions on curvature and distance has been considered by several authors including D. Yang \cite{DY}, P. Gianniotis \cite{panagiotis} and R. Hochard \cite{hochard}.

\begin{theorem}[{Special case of results in \cite{ST2}.}]
Suppose $(M^3,g_0)$ is a complete Riemannian manifold, and $x_0\in M$ so that
\beql{conditions2}
\left\{
\begin{aligned}
\Ric_{g_0}\geq -\al_0\\
\VolB_{g_0}(x_0,1)\geq v_0>0.
\end{aligned}
\right.
\eeq
Then for all $r\geq 2$ there exists a Ricci flow $g(t)$ (incomplete) defined on
$B_{g_0}(x_0,r)$ for $t\in [0,T]$, with $g(0)$ equal to (the restriction of) $g_0$
such that 
\beql{stable_balls}
B_{g(t)}(x_0,r-1)\subset\subset B_{g_0}(x_0,r),
\eeq
and with Ricci and volume bounds persisting in the sense that
\beql{newRicVol}
\left\{
\begin{aligned}
\Ric_{g(t)}\geq -\al\\
\VolB_{g(t)}(x_0,1)\geq v>0,
\end{aligned}
\right.
\eeq
but additionally so that 
\beql{curv_decay_conc}
|\Rm|_{g(t)}\leq \frac{c_0}{t},
\eeq
where $T, \al, v, c_0>0$ depend only on $\al_0$, $v_0$ and $r$.
Moreover, we have distance function convergence $d_{g(t)}\to d_{g_0}$ as $t\downto 0$
in the sense that (for example) for all $x,y\in B_{g_0}(x_0,r/8)$ we have
\beql{distest1}
d_{g_0}(x,y)-\be\sqrt{c_0 t}\leq
d_{g(t)}(x,y)\leq
e^{\al t} d_{g_0}(x,y)
\eeq
and
\beql{distest2}
d_{g_0}(x,y)\leq \ga [d_{g(t)}(x,y)]^{\frac{1}{1+4c_0}},
\eeq
where $\be$ is universal, $\al$ is as above, and $\ga$ depends only on $c_0$, thus on $\al_0$, $v_0$ and $r$.
\label{locRF}
\end{theorem}
Example \ref{ex3} illustrates the importance of \eqref{stable_balls} in the context of incomplete Ricci flows. Our flows have stable balls that do not instantly collapse.
This property, together with the refined information from \eqref{distest1} and \eqref{distest2} 
relies crucially on the  curvature estimates of \eqref{newRicVol} and \eqref{curv_decay_conc}.
In fact, it will be useful to split off what control one obtains, in slightly greater generality, in the following remark.

\brmk[Distance estimates.]
Given a Ricci flow $g(t)$ of arbitrary dimension $n$, 
defined only for \emph{positive} times  $t\in (0,T]$, 
and satisfying the curvature bounds
\begin{equation}
\left\{
\begin{aligned}
& \Ric_{g(t)}\geq -\al 
\\
& \Ric_{g(t)}\leq \frac{(n-1)c_0}{t}
\end{aligned}
\right.
\end{equation}
we automatically obtain curvature estimates like in Theorem \ref{locRF}
that are strong enough to extend the distance metric to $t=0$. 
Indeed, for any $0< t_1\leq t_2\leq T$, we have
\beq
\label{main_dist_est2}
d_{g(t_1)}(x,y)-\be\sqrt{c_0}(\sqrt{t_2}-\sqrt{t_1})
\leq d_{g(t_2)}(x,y)
\leq e^{\al(t_2-t_1)}d_{g(t_1)}(x,y),
\eeq
where $\be=\be(n)$.
In particular, the distance metrics $d_{g(t)}$ converge locally uniformly to some distance metric $d_0$ 
as $t\downto 0$, and 
\beq
\label{specialised_main_dist_est2}
d_0(x,y)-\be\sqrt{c_0t}
\leq d_{g(t)}(x,y)
\leq e^{\al t}d_0(x,y),
\eeq
for all $t\in (0,T]$.
Furthermore, there exists $\eta>0$ depending on $n$, $c_0$ and 
the size of the region in space-time we are considering
such that
\beq
\label{holder_est}
d_{g(t)}(x,y)\geq \eta \left[d_0(x,y)\right]^{1+2(n-1)c_0},
\eeq
for all $t\in (0,T]$.
Here we are being vague as to where the points $x,y$ live, and the dependencies of $\eta$.
In the local version of this result that we need, these issues require some care. See \cite[Lemma 3.1]{ST2} for full details.
\label{dist_est_rmk}
\ermk


\brmk
The distance estimates of Theorem \ref{locRF} imply that over appropriate local regions $\Om$, we have
$id:(\Om,d_{g(t)})\to (\Om,d_{g_0})$ is H\"older and 
$id:(\Om,d_{g_0})\to (\Om,d_{g(t)})$ is Lipschitz.
Similar statements follow from the estimates of  Remark \ref{dist_est_rmk}, with $d_0$ in place of $d_{g_0}$.
In particular, the identity is a bi-H\"older homeomorphism from the potentially-rough initial metric to the positive-time smooth metric.
\ermk

\brmk[Parabolic rescaling of Ricci flows.]
It will be essential to digest the natural parabolic rescaling of a Ricci flow to obtain a new Ricci flow. 
Given a Ricci flow $g(t)$ and a constant $\la>0$, we can define a scaled Ricci flow
$$g_\la(t):=\la g(t/\la),$$
on the appropriately scaled time interval.
When we do this rescaling, the curvature is scaled, for example the
sectional curvatures are all multiplied by a factor $\la^{-1}$.
One significant aspect of the uniform curvature estimate
\eqref{curv_decay_conc} obtained in Theorem \ref{locRF} is that it is \emph{invariant} under rescaling. That is, we get the same estimate for the rescaled flow \emph{with the same $c_0$}.
\label{parabolic_rescaling}
\ermk

\brmk
\label{higher_reg_rmk}
A well-known consequence of the curvature decay \eqref{curv_decay_conc},
via parabolic regularity, known as Shi's estimates in this context
\cite[\S 13]{formations}, is that for all $k\in \N$ we have
\beql{curv_decay_conc2}
|\grad^k\Rm|_{g(t)}\leq \frac{C}{t^{1+k/2}},
\eeq
where $C$ is allowed to depend on $k$ as well as $\al_0$, $v_0$ and $r$.
\ermk

\cmt{purposefully brushing under the carpet that this is an interior estimate, because it is just a distraction.}

\subsection{Getting the flow going}
\label{getflowgoing}

\cmt{recap for lecture 4:
Recall that we would like to start the flow with a RLS in 3D (locally, e.g. on ball of radius $r$). i.e. find RF $(M,g(t))$ for $t\in (0,T)$ so that $d_{g(t)}\to d_0$ as $t\downto 0$ in strong enough sense, with
$(M,d_0)$ isometric to the RLS.
If we can get good enough attainment of the initial data (cf. distance estimates we had last time) then 
over appropriate local regions $\Om\subset M$, we have
$id:(\Om,d_{g(t)})\to (\Om,d_0)$ is H\"older and 
$id:(\Om,d_0)\to (\Om,d_{g(t)})$ is Lipschitz.
This would allow us to show that the RLS is (locally) bi-H\"older
to $(\Om,d_{g(t)})$, i.e. $(\Om,{g(t)})$, via the identity map.
We sketched last time how all this would follow from proving existence of (and estimates for) a local Ricci flow (put more in notes?).
Thus today we focus on the local RF theorem. RESTATE IT}

The target of Theorem \ref{locRF} is to start the Ricci flow locally for a controlled amount of time, with a selection of quantitative estimates. But first we want to focus on the problem of merely starting the flow locally, without estimates and for an uncontrolled time interval. 
We follow the general strategy we used in the 2D theory, as in Section \ref{2Dsect}. To flow on a ball $B_{g_0}(x_0,r)$ as in the theorem, we would like to try to conformally blow up the metric in order to obtain a complete metric with constant negative curvature asymptotically, and hence one of bounded curvature. In 2D, this construction is straightforward \cite{JEMS}. In higher dimensions we use a lemma of Hochard \cite{hochard}, developed precisely for this purpose, following Simon \cite[Theorem 8.4]{simon2012}.

\cmt{If short of time, use only following paragraph, not the lemma. But we NEED the factor $\ga$ in the proof later}

The following lemma basically says that given a local region of a Riemmanian manifold with curvature bounded by $|\Rm|\leq \rho^{-2}$
on that local region, we can conformally blow up the metric near the boundary so that the resulting metric is \emph{complete} and has curvature bounded by $|\Rm|\leq \ga\rho^{-2}$, where the factor $\ga$ only depends on the dimension. It is asymptotically like  a scaled hyperbolic metric.

\begin{lemma}[{Variant of Hochard,  \cite[Lemma 6.2]{hochard}.}]
Let $(N^n,g)$ be a smooth Riemannian manifold (that need not be complete) and let $U \subset N$ be an open set. 
Suppose that for some $\rho\in (0, 1]$, we have
$\sup_{ U } |\Rm|_g  \leq \rho^{-2}$,
$B_{g}(x,\rho ) \subset \subset N$ and 
$\inj_{g}(x) \geq \rho$ for all $x \in U$. 
Then  there exist  
an open set $\ti U \subset  U$,
a smooth metric $\ti g$ defined on $\ti U$ and
a constant $\ga=\ga(n)\geq 1$,
such that 
each connected component of
$(\ti U, \ti g)$ is a complete Riemmanian  manifold such that
\begin{enumerate}[\ (1)]
\item
$\sup_{\ti U} |\Rm|_{\ti g} \leq  \ga\rho^{-2}$,
\item
$U_{\rho} \subset \ti U \subset  U$, and 
\item
$\ti g= g  \ \mbox{\rm on} \  \ti U_{\rho} \supset U_{2\rho}$,
\end{enumerate} 
where  
$U_s = \{ x \in U \ | \ B_{g}(x,s) \subset \subset U \}$.
\label{HochardsLemma2}
\end{lemma}
In order to achieve our modest task of at least starting the flow locally on a ball 
$B_{g_0}(x_0,r)$, we can appeal to Lemma \ref{HochardsLemma2} with $(N,g)$ of the lemma equal to 
$(M,g_0)$, with $U$ of the lemma equal to $B_{g_0}(x_0,r+2)\subset M$, and with $\rho\in (0,1]$
sufficiently small so that 
$\sup_{ U } |\Rm|_g  \leq \rho^{-2}$ and 
$\inj_{g}(x) \geq \rho$ for all $x \in U$. 
The output of the lemma is a complete Riemannian manifold $(\ti U,\ti g)$ of bounded curvature, where
$\ti U\subset M$, such that $\ti g=g_0$ on $B_{g_0}(x_0,r)$.
To obtain the local Ricci flow, we simply apply Shi's theorem \ref{shi_thm} to 
$(\ti U,\ti g)$ to obtain a Ricci flow $g(t)$, $t\in [0,T]$, and then restrict the flow
to $B_{g_0}(x_0,r)$.

We have succeeded in starting the flow, but there is a major problem that the existence time $T$ depends on $\sup_{\ti U}|\Rm|_{\ti g}$, which in turn depends on 
$\sup_{B_{g_0}(x_0,r+2)}|\Rm|_{g_0}$. This is not permitted. To be useful, the existence time can only ultimately depend on a lower Ricci bound, a lower bound on the volume of 
$B_{g_0}(x_0,1)$ and the radius $r$.

To resolve this problem, we need to extend the flow we have found. This will be acheived 
with the extension lemma of the following section, after which we'll pick up the proof of Theorem \ref{locRF} in Section \ref{use_of_ext_lem_sect}.


\subsection{The Extension Lemma}

The extension lemma will tell us that if we have a reasonably controlled Ricci flow on some local region $B_{g_0}(x_0,r_1)$, with reasonably controlled initial metric, then we can extend the Ricci flow to a longer time interval, albeit on a smaller ball, whilst retaining essentially the same estimates on the extension. This will be iterated to prove Theorem \ref{locRF}.

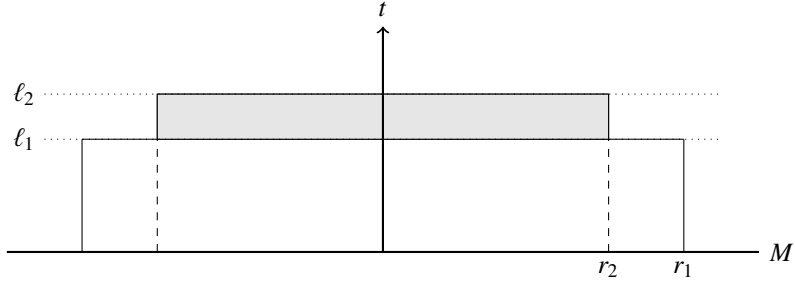
\begin{figure}[h]
\centering
\begin{tikzpicture}  

\newcommand\tscale{0.7}
\newcommand\axisheight{3}


\draw (-4,0) -- (-4,\axisheight/2);
\draw (4,0) node[below]{$r_1$} -- (4,\axisheight/2);
\draw (-4,\axisheight/2) -- (4,\axisheight/2);
\draw [dashed, thin] (3,0) node[below]{$r_2$} -- (3,\axisheight/2);
\draw (3,\axisheight/2) -- (3,\tscale*\axisheight);
\draw [dashed, thin] (-3,0) -- (-3,\axisheight/2);
\draw (-3,\axisheight/2) -- (-3,\tscale*\axisheight);
\draw (-3,\tscale*\axisheight) -- (3,\tscale*\axisheight);


\draw[fill=black!10] (-3,\axisheight/2) rectangle (3,\tscale*\axisheight);


\draw [dotted, thin] (-4.5,\axisheight/2) node[left]{$\ell_1$} -- (4.5,\axisheight/2);
\draw [dotted, thin] (-4.5,\tscale*\axisheight) node[left]{$\ell_2$} -- (4.5,\tscale*\axisheight);

\draw [thick, <-] (0,\axisheight) node[above]{$t$} -- (0,0) -- (5,0) node[right]{$M$};
\draw [thick, -] (0,0) -- (-5,0);

\end{tikzpicture}
\caption{Domain of original Ricci flow In the Extension Lemma \ref{ext_lemma} and (shaded) additional domain to which we extend.}
\label{fig_ext_lem}
\end{figure}

\begin{lemma}[{Extension lemma \cite{ST2}.}]
Given $v_0>0$, there exist  $c_0\geq 1$ and $\tau>0$
such that the following is true.
Let $r_1 \geq 2$, and let $(M,g_0)$ be a smooth three-dimensional Riemannian manifold with  $B_{g_0}(x_0,r_1) \subset \subset M$, and
\begin{enumerate}[(i)]
\item
$\Ric_{g_0} \geq -\al_0$ for some $\al_0  \geq 1$ on $B_{g_0}(x_0,r_1)$, and  
\item
$\VolB_{g_0}(x,s) \geq v_0 s^3$ for all $s\leq 1$ and all $x \in B_{g_0}(x_0,r_1-s).$
\end{enumerate}
Assume additionally that we are given a smooth Ricci flow  $(B_{g_0}(x_0,r_1),g(t))$, $t \in [0,\ell_1]$, where 
$\ell_1 \leq \frac{\tau}{200\al_0 c_0}$,
with $g(0)$ equal to the restriction of $g_0$, for which
\begin{enumerate}[(a)]
\item
$ |\Rm|_{g(t)} \leq \frac{c_0}{t}$  and
\item
$\Ric_{g(t)}\geq -\frac{\tau}{\ell_1}$
\end{enumerate}
on $B_{g_0}(x_0,r_1)$  for all $t \in (0,\ell_1]$.
Then, setting $\ell_2  = \ell_1(1+\frac{1}{4c_0})$ and 
$r_{2} = r_{1} - 6\sqrt{\frac{\ell_{2}}{\tau}}\geq 1$,
the Ricci flow $g(t)$ can be extended smoothly to a Ricci flow over
the longer time interval $t\in [0,\ell_2]$, albeit on
the smaller ball $B_{g_0}(x_0, r_{2})$,  with 
\begin{enumerate}[(a$'$)]
\item
$ |\Rm|_{g(t)} \leq \frac{c_0}{t}$  and
\item
$\Ric_{g(t)}\geq -\frac{\tau}{\ell_2}$
\end{enumerate}
throughout $B_{g_0}(x_0,r_2)$  for all $t \in (0,\ell_2]$.
\label{ext_lemma}
\end{lemma}
See Figure \ref{fig_ext_lem} for the domains where the Ricci flows are defined.


\subsection{Using the Extension lemma to construct a local Ricci flow}
\label{use_of_ext_lem_sect}

Assuming the Extension lemma, we can now pick up from where Section \ref{getflowgoing} 
ended,  and give a sketch of the proof of 
the existence and curvature bounds of Theorem \ref{locRF}.
See the end of Section 4 of \cite{ST2} for the fine details.
For the estimates on the distance function, the crucial thing is to have the lower Ricci and upper $\Rm$ bounds of the theorem. The estimates then follow from \cite[Lemma 3.1]{ST2}, cf. Remark \ref{dist_est_rmk}.

Before we start, it will be useful to reformulate the conditions \eqref{conditions2} of Theorem \ref{locRF} a little.
The claim is that taking $r$ from Theorem \ref{locRF}, the conditions \eqref{conditions2} imply 
that
for $r_1=r+1$, 
\beql{implied_cond}
\VolB_{g_0}(x,s) \geq \ti v_0 s^3\quad\text{  for all }s\leq 1\text{ and all }x \in B_{g_0}(x_0,r_1-s),
\eeq
where $\ti v_0>0$ depends on $v_0$, $\al_0$ and  $r$. 
To prove this is a simple comparison geometry argument; if we had 
\eqref{conditions2} with $\al_0=0$, then this was done in 
Section \ref{comp_geom_sect} -- see \eqref{vol_rat_lower_bd} in particular. The general case is a minor variant.

\cmt{GO THROUGH ARGUMENT WHEN SPACE IS RESTRICTED!!!}
\cmt{for the argument above, work on the ball of radius $r_1+1$. We need some space! Or I guess we could mimic the argument from before and work on a ball of much larger radius.}

This observation means that we may as well replace the second condition
of \eqref{conditions2} by the condition \eqref{implied_cond}. We will do that below, dropping the tilde on $\ti v_0$.
Of course, this means that we satisfy 
the conditions $(i)$ and $(ii)$ of Lemma \ref{ext_lemma}.

Now we have this new $v_0$, we can use it to obtain constants $c_0$ and $\tau$ from Extension lemma \ref{ext_lemma}.

The next step is to get the flow going as in Section \ref{getflowgoing}, except that we construct the flow on the ball $B_{g_0}(x_0,r+1)$ of slightly larger radius than before. 
This flow will exist on some nontrivial but otherwise uncontrolled time interval $[0,\ell_1]$. By reducing $\ell_1>0$ if necessary we can ensure that 
conditions $(a)$ and $(b)$ of the Extension lemma are satisfied.

At this point, if $\ell_1 > \frac{\tau}{200\al_0 c_0}$, then we have existence of a Ricci flow on a time interval with the right dependencies. But in general, we will have existence for an uncontrollably short time interval. In this latter case, we apply the Extension lemma.
This gives us existence over a longer time interval, albeit on a smaller ball, and with the same estimates. Thus we can iterate this process to extend the flow more and more. 

We stop iterating only if either $\ell_1$ is large enough as above, or if the radius of the ball on which we get estimates becomes smaller than the radius $r$ on which we would like our ultimate solution. Here is a potential worry: If the original existence time were extremely small, then we would need a huge number of iterations to get a reasonable existence time. Each time we extend on a smaller ball, so there is a concern that the whole space is quickly exhausted as the balls get too small. Luckily, the smaller the original existence time, the smaller the reduction in the size of the ball need be. It turns out that we can ask that the radius of the balls where we obtain existence never drops below $r$, but still be guaranteed a good time of existence.

Section 4 of \cite{ST2} gives the details.


\subsection{{Proof of the Extension lemma \ref{ext_lemma}}}

We sketch some of the main ingredients and ideas of the proof of the Extension lemma.
The objective is to make the full proof in \cite[\S 4]{ST2} accessible.



The first result expresses that a Ricci flow that has an initial volume lower bound, 
and a Ricci lower bound for each time, must regularise itself over short time periods.
The lemma has many variants; in particular, a global version can be found in
the earlier work of Miles Simon \cite{simon2012}.

\begin{lemma}[{Local lemma, cf. \cite[Lemma 2.1]{ST1}.}]
Let $(N^3,g(t))$, for ${t\in [0,T]}$, be a  Ricci flow such that for some fixed $x\in N$ we have $B_{g(t)}(x,1)\subset\subset N$ for all $t\in [0,T]$, and so that
\begin{enumerate}[(i)]
\item
$\VolB_{g(0)}(x,1) \geq v_0>0 $, and 
\item
$\Ric_{g(t)} \geq -1$ on $ B_{g(t)}(x,1) $ for all $t\in [0,T]$. 
\end{enumerate}
Then there exist $C_0 \geq 1 $ and $\hat T >0$, depending only on $v_0$, such that 
$ |\Rm|_{g(t)}(x) \leq 
C_0/t
$,  and
$\inj_{g(t)}(x) \geq \sqrt{
t/C_0
}$ 
for all $0<t \leq \min(\hat T,T).$
\label{modifiedlocal2}
\end{lemma}
The key point to digest here is that we obtain the curvature decay 
$|\Rm|_{g(t)} \leq C_0/t$, which is an estimate that is invariant under parabolic rescaling of the flow
as in Remark \ref{parabolic_rescaling}. 
The main idea of the proof is as follows. If the curvature estimate failed in a shorter and shorter time and with disproportionately larger and larger $C_0$, for some sequence of flows, then we can 
parabolically rescale each flow so that they exist for a longer and longer time, now with almost nonnegative Ricci curvature, but still end up with $|\Rm|=1$ at the end.
It then makes sense to pass to a limit of the flows (in the smooth Cheeger-Gromov sense, see, for example, \cite[\S 7]{RFnotes} and the discussion in Section \ref{how_to_apply_sect})
and obtain a limit smooth Ricci flow that existed since time $-\infty$, has $\Ric\geq 0$,
and has $|\Rm|=1$ at time $0$, say, but which has bounded curvature.
Moreover it is possible to argue that the hypothesis $\VolB_{g(0)}(x,1) \geq v_0>0 $
eventually tells us that this limit Ricci flow has positive asymptotic volume ratio in the sense
that the volume of a ball of radius $r>0$ is at least $\ep r^3$ for some $\ep>0$, however large we take $r$.
Finally, a result of Chow and Knopf \cite[Corollary 9.8]{ChowKnopf} tells us that such a 3D complete ancient Ricci flow with bounded curvature must have nonnegative sectional curvature.

This puts our limit flow in the framework of Perelman's famous $\kappa$-solutions \cite[\S 11]{P1} that are pivotal in understanding the structure of Ricci flow singularities. 
The output of Perelman's work is that the asymptotic volume ratio of such a nonflat solution must in fact be \emph{zero}, contradicting what we already saw, and completing the (sketch) proof.
See \cite[\S 5]{ST1} for the argument in full detail.


The following key ingredient gives a local version of the preservation of nonnegative Ricci curvature that we saw at the beginning of these lectures. The Ricci lower bounds it provides will be critical in the applications we have in mind.

\begin{lemma}[{cf. Double Bootstrap Lemma \cite[Lemma 9.1]{ST1}.}]
Let $(N^3,g(t))$, for ${t\in [0,T]}$, be a Ricci flow such that 
$B_{g(0)}(x,2)\subset\subset N$ for some $x\in N$,
and so that throughout $B_{g(0)}(x,2) $ we have
\begin{enumerate}[(i)]
\item
$|\Rm|_{g(t)} \leq \frac{c_0}{t}$  for some $c_0\geq 1$ and all $t\in (0,T]$, and 
\item
$\Ric_{g(0)} \geq -\de_0$ for some $\de_0>0$.
\end{enumerate}
Then there exists $\hat S = \hat S(c_0,\de_0)>0$ such that 
for all $0\leq t \leq \min(\hat S,T)$ we have
$$ \Ric_{g(t)}(x) \geq - 100 \de_0 c_0. $$ 
\label{modifieddouble2}
\end{lemma}
Equipped with these two lemmata, we will attempt to outline some of the central ideas of the Extension lemma \ref{ext_lemma}. Amongst the various details we suppress is the estimation of the exact domain on which we can extend the flow. For the precise details of the proof, see Section 4 of \cite{ST2}.

\vskip 0.2cm
\noindent
{\bf Step 1: Choose the constants}

Let $v_0$ be the constant from the Extension lemma \ref{ext_lemma},
and plug it in to the 
Local lemma \ref{modifiedlocal2} to obtain constants $C_0$ and $\hat T$.
Extracting the constant $\ga\geq 1$ from Hochard's lemma \ref{HochardsLemma2}, we are already in a position to specify $c_0:=4\ga C_0$.

In turn this allows us to specify $\de_0:=\frac{1}{100c_0}$, and appeal to the Double Bootstrap lemma \ref{modifieddouble2} to obtain $\hat S$. We can then insist to work on time intervals 
$[0,\tau]$ where $\tau\leq \hat T, \hat S$ so that all our supporting lemmata apply.

\vskip 0.2cm
\noindent
{\bf Step 2: Improvement of the curvature decay}

The flow we are given in the Extension lemma already enjoys $c_0/t$ curvature decay, where
we recall that $c_0:=4\ga C_0>C_0$. By parabolically scaling up our flow, scaling up time by a factor $\tau/\ell_1$, and distances by a factor $\sqrt{\tau/\ell_1}$, our flow ends up with a lower Ricci bound of $-1$, and the Local lemma \ref{modifiedlocal2} applies over all unit balls to give $C_0/t$ curvature bounds (and injectivity radius bounds) for $t\in (0,\tau]$. Scaling back, we get again $C_0/t$ curvature bounds for the original flow, holding now
up until time $\tau \times \frac{\ell_1}{\tau}=\ell_1$, as desired.

Hidden a little here is that the details of this step require the assumed $c_0/t$ curvature decay in order to verify the hypotheses of the Local lemma.

\vskip 0.2cm
\noindent
{\bf Step 3: Continuation of the flow}

The previous step has obtained improved control on the curvature of the originally given Ricci flow on an appropriate interior region of $B_{g_0}(x_0,r_1)$, and in particular at time $\ell_1$ we have
$$|\Rm|_{g(\ell_1)}\leq \frac{C_0}{\ell_1}.$$
We can now follow Hochard, Lemma \ref{HochardsLemma2}, and modify the metric around the outside in order to make it complete and with a curvature bound
$$|\Rm|\leq \ga \frac{C_0}{\ell_1}=\frac{c_0}{4\ell_1},$$
by definition of $c_0$.
Having a complete, bounded curvature metric, we can apply Shi's theorem \ref{shi_thm}
to extend the flow for a time $\frac{1}{16}\times \frac{4\ell_1}{c_0}=\frac{\ell_1}{4c_0}$,
maintaining a curvature bound 
$$|\Rm|\leq 2\times \frac{c_0}{4\ell_1}\leq \frac{c_0}{\ell_2},$$
where $\ell_2:=\ell_1(1+\frac{1}{4c_0})\leq 2\ell_1$. Thus we have proved a little more than 
$$|\Rm|_{g(t)}\leq \frac{c_0}{t}$$
for the extended Ricci flow $g(t)$, which is part $(a')$ of the desired conclusion.

\vskip 0.2cm
\noindent
{\bf Step 4: Ricci lower bounds}

Now we have indirectly obtained $c_0/t$ curvature decay for the extended flow, we are in a position to obtain the crucial lower Ricci bound conclusion of part $(b')$.
Similarly to before, we parabolically scale up the extended flow, scaling time by a factor 
$\tau/\ell_2$ so that the new flow lives for a time $\tau$. This preserves our $c_0/t$ curvature decay, but improves the initial Ricci lower bound from $-\al_0$ to 
$$-\al_0\frac{\ell_2}{\tau}\geq -2\al_0\frac{\ell_1}{\tau}
\geq -\frac{1}{100c_0}=-\de_0.$$
Thus we can apply the Double Bootstrap lemma \ref{modifieddouble2} to conclude that the rescaled flow has a Ricci lower bound of $-100c_0\de_0=-1$.
Rescaling back, we obtain that our original extended flow has a Ricci lower bound of 
$-\tau/\ell_2$ as required.

\cmt{Wordy summary:
We started with control on $[0,\ell_1]$. We improved the control using the local lemma. We did the conformal hyperbolic blow-up trick near the boundary, giving up some of that improved control. We flowed forwards with Shi. We analysed the extended flow using double bootstrap in order to end up with the same control as before, but on a longer time interval!
}

This completes the sketch of the proof of the Extension lemma \ref{ext_lemma}.

\brmk
The Pyramid Ricci flow construction that we briefly discuss in Section \ref{pyramid_sect}
requires its own extension lemma. In that result we are unable to follow this strategy of applying Shi's theorem to extend, since the curvature bounds available are insufficient. 
\ermk


\section{Pyramid Ricci flows}
\label{pyramid_sect}

In  Section \ref{localRFsect} we stated and sketched the proof of 
the Local Ricci flow theorem \ref{locRF}. This result is enough to prove Theorem \ref{ACCT3Dthm} and hence settle the three-dimensional Anderson-Cheeger-Colding-Tian conjecture, but instead of elaborating on the details, in this section we give a different Ricci flow existence result that will not only give the slightly stronger result Theorem \ref{pyramid_version}, but will allow us a slightly less technical description of the proof.

The existence result we will give is for so-called Pyramid Ricci flows. To understand the idea, let us revisit the Local Ricci flow theorem \ref{locRF} in the case that the radius $r$ is a natural number $k$. 
For each $k\in\N$ we obtain a Ricci flow on a parabolic cylinder
$B_{g_0}(x_0,k)\times [0,T_k]$, and in general the final existence time $T_k$ converges to zero as $k\to\infty$
as demonstrated by Example \ref{ex4}.
The Ricci flows corresponding to different values of $k$ will exist on cylinders with nontrivial intersection, but it is unreasonable to expect that they agree on these intersections. 

The point of Pyramid Ricci flows is to fix this problem. Instead of obtaining a countable number of Ricci flows on individual cylinders, 
indexed by $k$, we obtain one Ricci flow that lives throughout
the \emph{union} of all cylinders.

\begin{theorem}[{with McLeod \cite{MT}, using \cite{ST1, ST2}.}]
Suppose $(M^3,g_0)$ is a complete Riemannian manifold, and $x_0\in M$ so that
\beql{conditions3}
\left\{
\begin{aligned}
\Ric_{g_0}\geq -\al_0\\
\VolB_{g_0}(x_0,1)\geq v_0>0.
\end{aligned}
\right.
\eeq
Then there exist increasing sequences $C_j \geq 1$ and $\al_j > 0$ and a decreasing sequence $T_j > 0,$ all defined for $j \in\N$, and depending only on $\al_0$ and $v_0$, for which the following is true.
There exists a smooth Ricci flow $g(t)$, defined on a subset of spacetime 
$$\cd_k:=\bigcup_{k\in\N} B_{g_0}(x_0,k)\times [0,T_k],$$
satisfying that $g(0)=g_0$ throughout $M$, and that for each $k \in \N,$
\beql{pyr_conc} 
\left\{
\begin{aligned}
\Ric_{g(t)}\geq -\al_k\\
|\Rm|_{g(t)}\leq \frac{C_k}{t},
\end{aligned}
\right.
\eeq
throughout $\B_{g_0} ( x_0 , k) \times \left( 0 , T_k \right]$.
\label{pyramid_thm}
\end{theorem}
The idea of defining a Ricci flow on a subset of spacetime was used by Hochard \cite{hochard}, prior to this work, when he defined
\emph{partial} Ricci flows. Our flows differ in that we concede having a uniform $c/t$ curvature bound in return for strong uniform control on the shape of the spacetime region on which the flow is defined, which is essential for the compactness arguments that use this theorem.

The theory that allows this pyramid construction revolves around a variant of the Extension lemma \ref{ext_lemma} known as the Pyramid Extension Lemma; see \cite[Lemma 2.1]{MT} and \cite[Lemma 4.1]{MT2}.

\section{{Proof of the 3D Anderson-Cheeger-Colding-Tian conjecture}}
\label{ACCTsect}

In this section we apply what we have discovered about Ricci flow
in order to complete the programme 
of Section \ref{how_to_apply_sect} and start the Ricci flow with a Ricci limit space in order to prove Theorems \ref{ACCT3Dthm} and \ref{pyramid_version}.

The theory can be carried out with Theorem \ref{locRF}, but here we will describe the strategy of \cite{ST2} using instead Theorem \ref{pyramid_thm}, with one benefit being a slight reduction in the technical aspects of considering distance metrics on incomplete manifolds.

\cmt{decided to be vague about the CG compactness below. We've given earlier some references for the usual CG convergence and compactness}

To reiterate the argument, we use Theorem \ref{pyramid_thm} to flow each $(M_i,g_i)$ generating our Ricci limit space $(X,d)$.
We call the resulting pyramid Ricci flows $g_i(t)$.
The  estimates on the curvature from Theorem \ref{pyramid_thm}
allow us to invoke a local version of Cheeger-Gromov-Hamilton compactness to get a type of smooth convergence
\beql{sm_con_again}
(M_i,g_i(t),x_i)\to (M,g(t),x_\infty)
\eeq 
to some smooth limit Ricci flow $g(t)$, which inherits the curvature bounds of the pyramid flows $g_i(t)$ and lives on a subset of spacetime similar to the domains of the pyramid flows. By keeping track of the evolution of distances, see 
Remark \ref{dist_est_rmk}, 
and \cite[Lemma 3.1]{ST2} for the precise local statement,
we obtain the following.

\begin{theorem}[Ricci flow from a Ricci limit space.]
Given a Ricci limit space $(X,d,\hat x_\infty)$ that is approximated by 
a sequence $(M^3_i , g_i ,x_i )$ of complete, smooth, pointed Riemannian three-manifolds such that 
\beql{conditions4}
\left\{
\begin{aligned}
\Ric_{g_i}\geq -\al_0\\
\VolB_{g_i}(x_i,1)\geq v_0>0,
\end{aligned}
\right.
\eeq
there exist increasing sequences $C_k \geq 1$ and $\al_k > 0$ and a decreasing sequence $T_k > 0,$ all defined for $k \in\N$, and depending only on $\al_0$ and $v_0,$ for which the following holds.

There exist a smooth three-manifold $M$, a point $x_\infty \in M$, 
and a complete distance metric $d_0 : M \times M \rightarrow [0,\infty)$ generating the same topology as we already have on $M$,
such that $(M,d_0,x_\infty)$ is isometric to the Ricci limit space
$(X,d,\hat x_\infty)$. 

Moreover, there exists a smooth Ricci flow $g(t)$ 
defined on a subset of spacetime $M \times (0,\infty)$ that contains $\B_{d_0} ( x_0 , k) \times \left(0, T_k \right]$ for each $k \in \N,$ 
with $d_{g(t)} \rightarrow d_0$ locally uniformly on $M$ as $t \downarrow 0$.
%
For any $k \in \N$  we have
\begin{equation}
\label{r_d_concs} 
\left\{
\begin{aligned}
\Ric_{g(t)} \geq -\al_k\\
| \Rm |_{g(t)} \leq \frac{C_k}{t}
\end{aligned}
\right.
\end{equation}
throughout $B_{d_0} ( x_0 , k ) \times (0,T_k]$.
Finally, for all $\Om\subset\subset M$, there exists $\si>1$
depending on $\Om$, $\al_0$ and $v_0$ 
such that for $t\in (0,\min\{\frac{1}{\si},T\}]$ and 
for all $x,y\in \Om$,  we have
\beql{simp_dist}
\frac{1}{\si}d_{g(t)}(x,y)\leq
d_{g_0}(x,y)\leq
\si [d_{g(t)}(x,y)]^\frac{1}{\si}.
\eeq
\label{pyramid_limit_flow}
\end{theorem}
Of course, the existence of the distance metric $d_0$ and its relation to the original Ricci limit space is established by virtue of having the Ricci flow $g(t)$, not the other way round as the phrasing above may suggest.

The simplified distance estimate \eqref{simp_dist} tells us immediately that for sufficiently small $t\in (0,T]$, the identity map 
$id:(\Om,d_{g(t)})\to (\Om,d_{g_0})$ is H\"older and 
$id:(\Om,d_{g_0})\to (\Om,d_{g(t)})$ is Lipschitz.

This is a very rough outline. See \cite{ST2, MT} for further details,
and \cite{MT2} for a discussion of subsequent developments and simplifications.

\section{An open problem}

Despite all the progress in this general area, there are some very basic problems that remain open.
An example of an intriguing long-standing open problem at the heart of the matter, which contrasts with Example \ref{ex4}, is the following.

\begin{conjecture}
Given a complete Riemannian three-manifold with $\Ric\geq 0$, there exists a smooth complete Ricci flow continuation.
\end{conjecture}

\section*{Appendix - Gromov compactness}
\label{grom_appen}
%
%

When defining Ricci limit spaces we used the following elementary but important compactness result.

\begin{theorem}[Gromov compactness, special case]
\label{grom_cpt_thm}
Given a sequence of complete pointed Riemannian $n$-manifolds $(M_i,g_i,x_i)$ with a uniform lower Ricci curvature bound, we may pass to a subsequence to obtain pointed Gromov-Hausdorff convergence
$$(M^n_i,g_i,x_i)\to (X,d,x_\infty),$$
to some pointed limit metric space $(X,d,x_\infty)$.
\end{theorem}
To illustrate the (simple) proof, let's see how one can prove the even simpler statement that for any $r>\ep>0$, there exists a metric space $(X,d,x_\infty)$ such that after passing to a subsequence, all elements of our sequence $(M^n_i,g_i,x_i)$ are $\ep$-close to 
$(X,d,x_\infty)$ in the sense that there exist maps $f_i:B_{g_i}(x_i,r)\to X$ such that 
\begin{enumerate}
\item
$f_i(x_i)=x_\infty$,
\item
for all $x,y\in B_{d_i}(x_i,r)$ we have
$$|d(f_i(x),f_i(y))-d_i(x,y)|<\ep,$$
\item
the $\ep$-neighbourhood of the image $f_i(B_{d_i}(x_i,r))$
contains the ball $B_d(x_\infty, r-\ep)\subset X$.
\end{enumerate}
To see this, for each fixed $i\in \N$ pick a maximal set of points
$x_i=:p^1_i, p^2_i,\ldots, p^m_i$ in $B_{g_i}(x_i,r)$ such that the balls 
$B_{g_i}(p^j_i,\ep/9)$ 
lie in $B_{g_i}(x_i,r)$ and are pairwise disjoint. 
Maximal means that we cannot adjust the points to squeeze in even one more point with this disjointness property.
That way, we can be sure that the balls $B_{g_i}(p^j_i,\ep/3)$ cover all
of $B_{g_i}(x_i,r)$.

\cmt{we've gone from $\ep/9$ to $\ep/3$ to handle the usual boundary issues. If we had a point not in the $\ep/3$ supposed-covering, and it was near the boundary, then we have to move it in towards the centre by $\ep/9$ and then take the $\ep/9$ around it.}

The claim is that the number of such points is necessarily bounded independently of $i$. Assuming this is true for the moment, we can pass to a subsequence so that for each $i$ we have the same number $m$ of these points. The distance between each fixed pair of points is bounded 
above by $2r$, and below by $2\ep/9$, so by passing to enough subsequences (but finitely many) we have convergence of the distances
$$d_i(p^j_i,p^k_i)\to d^{jk}>0$$
as $i\to\infty$, where $j\neq k$.
The metric space $X$ can be chosen to be a finite collection of points
$X=\{q^1,\ldots, q^m\}$, with distance
$$d(q^j,q^k):=d^{jk},$$
and we can choose $x_\infty:=q^1$ and a suitable map $f_i$ such that $f_i(p^j_i)=q^j$, giving
$$|d(f_i(p_i^j),f_i(p_i^k))-d_i(p_i^j,p_i^k)|=|d(q^j,q^k)-d_i(p_i^j,p_i^k)|\to 0,$$
as required.

\cmt{shame that taking the GH approximations the other way would clash with the definition we took}

To show that the number of points $\{p^j_i\}$ is bounded independently of $i$, we use Bishop-Gromov and consider volumes.
On the one hand, the lower Ricci bound and volume comparison (cf. Theorem \ref{BGthm} and Remark \ref{gen_al0_rmk}) tell us that the volume of $B_{g_i}(x_i,r)$ is uniformly bounded above independently of $i$. 
On the other hand, an argument similar to that giving \eqref{vol_rat_lower_bd} (but extended to handle
lower Ricci bounds other than zero) tells us that the volume of each pairwise disjoint ball
$B_{g_i}(p^j_i,\ep/9)$ is uniformly bounded below by some $i$-independent positive number. 
Combining these two facts gives an $i$-independent upper bound on the total number of such balls that can be pairwise disjoint and fit within $B_{g_i}(x_i,r)$.

For further aspects of the proof, see for example \cite[Chapter 3]{cheeger_book}.

\parskip 1pt

\setlength{\itemsep}{0pt}

\end{document}